\documentclass[11pt]{article}
\usepackage{amssymb,amsmath,amstext,amscd,amsthm,
graphicx,fancyhdr,rotating}
\theoremstyle{plain}
\newtheorem{lemma}{Lemma}
\newtheorem{defn}{Definition}
\newtheorem{prop}[lemma]{Proposition}
\newtheorem{theo}{Theorem}
\newtheorem{cor}{Corollary}
\newtheorem{remark}{Remark}
\newtheorem{exmp}{Example}

\begin{document}

\begin{titlepage}

\begin{flushright}
ICMPA-MPA/revf/2009/16\\
NITheP-09-10\\
\end{flushright}

\begin{center}

{\Large\bf Multi-indicial symmetric functions}

Joseph Ben Geloun$^{a,b,c,*}$ and Mahouton Norbert Hounkonnou$^{b,\dag}$

$^{a}${\em National Institute for Theoretical Physics (NITheP)}\\
{\em Private Bag X1, Matieland 7602, South Africa}\\
$^{b}${\em International Chair of Mathematical Physics
and Applications}\\
{\em ICMPA--UNESCO Chair, 072 B.P. 50  Cotonou, Republic of Benin}\\
$^{c}${\em D\'epartement de Math\'ematiques et Informatique}\\
{\em  Facult\'e des Sciences et Techniques, Universit\'e Cheikh Anta Diop, Senegal}

E-mails:  $^{*}$bengeloun@sun.ac.za,\quad $^{\dag}$norbert$_-$hounkonnou@cimpa.net

\begin{abstract}
In this paper, using the theory of category,
we generalize known properties of symmetric polynomials
and functions and characterize the multi-indicial symmetric functions.
Examples have been given on Schur functions.
\end{abstract}

\today
\end{center}


MCSs: 05E05, 13B25.

\end{titlepage}

\section{Introduction}

A great deal of attention has been paid to the symmetric functions
and orthogonal polynomials (\cite{IGM,GR,GR2} and references
therein). Indeed, symmetry is an inescapable feature of most physical
phenomena. Following \cite{IGM}, the theory of symmetric functions is one of
the most classical parts of algebra, going back to the $16^{th}$ and $17^{th}$
centuries and attempts of mathematicians of that epoch to solve polynomial
equations of degree higher than two. Generalization of symmetric functions
in several sets of variables (the so called multisymmetric functions)
was found by McMahon in the beginning of the past century \cite{mac}.
Still recently, McMahon symmetric polynomials have
been studied in different contexts \cite{gessel}-\cite{fresc}.
For instance in \cite{gessel}, the McMahon symmetric polynomials in two
sets of variables have been used to find explicit formulas and to prove
$P$-recursiveness for some objects such as Latin rectangles and $0-1$ matrices
with zeros on the diagonal and given row and column sums.
Thereafter, using the approach by McDonald \cite{IGM}, Dalbec extended the theory
of multisymmetric functions in two sets of variables to the multihomogeneous case,
the so called factorizable forms, in characteristic $0$ field and provided with
a MAPLE code for generating such objects \cite{dalbec}.
Vaccarino \cite{vaccar} generalized the above results  as well as those of \cite{fresc}
(dealing with characteristic $2$ fields) to the ring of multisymmetric
functions over a commutative ring.

Among the various families of symmetric functions,
the most significant are undoubtedly the Schur functions, because of their
intimate relationship with the irreducible characters of both the symmetric
group and the general linear groups, and for their combinatorial applications.

In this paper, the McDonald formalism has been extended using the
theory of category, in order to define multi-indicial symmetric functions
including different sets of variables with several tensorial indices.
More specifically, this paper addresses results on
two remarkable classes of symmetric functions with mixed types of tensor indices
and introduces their full characterization. Illustration has been given on Schur functions.

In Section 2, we give a generalization of known properties
of the ring of symmetric polynomials.
The ring of symmetric functions $\Lambda$ which is an inverse limit
is defined as a universal object. In Section 3,
we deal with the study of multi-indicial symmetric polynomials.
Relevant properties of the graded rings of such polynomials are derived.
The multi-indicial symmetric functions are logically introduced.
Section 4 is devoted to the definition of multi-indicial partition
and the corresponding definition of the Schur function.
We end the paper with some concluding remarks.

\section{Symmetric polynomials: main results}

In this section, we build the theoretical framework of our study. For that,
we recall main properties of the ring of symmetric polynomials
and give their generalization.
The ring of symmetric functions is defined as a universal object.

Let us introduce the definition \cite{IGM}:
\begin{defn} \label{ring}
Let ${\bf x}_1,{\bf x}_2,\ldots,{\bf x}_n$ be $n$  independent indeterminates,
$S_{n}$ be the symmetric group of permutations of a set with $n$ elements
acting on the polynomial ring ${\mathbb{Z}}\left[ {\bf x}_1, {\bf x}_2, \ldots ,
{\bf x}_n\right]$
by permuting the indeterminates, i.e:
\begin{eqnarray}
&&\forall P = a_{I}{\bf x}^{I} \in {\mathbb{Z}}\left[ {\bf x}_1, {\bf x}_2, \ldots ,
{\bf x}_n\right],\\
&&\forall \sigma \in S_n,\;\;\; \sigma P = \sigma a_{I}{\bf x}^{I} =
  a_{I}{\bf x}^{I}_{\sigma(\cdot)},
  \nonumber
\end{eqnarray}
where ${\bf x}={\bf x}_{1}{\bf x}_{2}\dots{\bf x}_{n}$, $a_{I} \in \mathbb{Z}$.
 $I= (i_1, i_2, \ldots, i_k)$,
 with $ 0\leq i_k$ and $\; 1 \leq k \leq n$,
denotes the usual multi-index notation (the implicit summation is used).
Then,
${\Lambda}_{n} :={\mathbb{Z}}\left[ {\bf x}_1, {\bf x}_2, \ldots ,{\bf x}_n \right]^{S_n}$
is the subring of $ {\mathbb{Z}}\left[ {\bf x}_1, {\bf x}_2, \ldots
,{\bf x}_n\right]$ of symmetric polynomials obtained by permuting the ${\bf x} _i$.
\end{defn}
\begin{remark} Let us pay attention to the fact that this sum is
globally invariant under any permutation, instead of the monomial terms taken
separately. For example, ${\bf x} _{J}$ may not be equal to ${\bf x} _{\sigma J}$.
\end{remark}
\begin{exmp}
Assume n=4, i.e the set of indeterminates is
$\left\{ {\bf x}_1, {\bf x}_2, {\bf x}_3, {\bf x}_4 \right\}$.
 The following polynomials belong to
$\Lambda_{4}$: $f^{1} = {\bf x}_1 + {\bf x}_2 + {\bf x}_3 + {\bf x}_4, \;\;\;
\forall r\in
\mathbb{N},\; f^{r} = {\bf x}_1 ^r + {\bf x} _2 ^r+ {\bf x} _3 ^r+ {\bf x}_4 ^r,
\;\;\; f =
{\bf x}_1{\bf x} _2 + {\bf x} _1{\bf x} _3 +{\bf x} _1{\bf x} _4 +
{\bf x} _2{\bf x} _3
+ {\bf x} _2{\bf x}_4 + {\bf x}_3{\bf x}_4.$
\end{exmp}
 If $f \;\in\; \Lambda _{n}$, one can write
$\;\;f = \sum_{r\geq 0} f^{r},$
where $f^{r}$ is the homogeneous component of $f$ of degree $r$.
One can verify that each of the $f^{r}$ is itself invariant under $S_n$ and
hence, $\Lambda _n$ is a graded ring. This statement can be written as:
$\;\;
\Lambda_{n} = \bigoplus _{r\geq 0} \Lambda_{n}^{r},
$
where $\Lambda_{n}^{r}$
is the additive group of homogeneous symmetric polynomials in
$\left\{ {\bf x} _1, {\bf x} _2,\ldots , {\bf x} _n\right\}$, provided the
following convention:
$0$ is homogeneous of any degree.  One requires also that
a polynomial of degree $0$ is nothing but an element of the coefficient
ring, i.e $\Lambda_{n}^0 = \mathbb{Z}$.

Adding a new indeterminate ${\bf x} _{n+1}$, we
can realize the ring
\[ \Lambda _{n+1} =
\mathbb{Z} \left[ {\bf x} _1, {\bf x} _2, \ldots,
 {\bf x}_{n}, {\bf x}_{n+1}\right] ^{S_{n+1}}\]
and the following statement holds \cite{IGM}.
\begin{lemma}\label{sim}
Let $\pi _{n+1}$ be the mapping from $\Lambda _{n+1}$ to $\Lambda _{n}$ defined
by setting ${\bf x} _{n+1} = 0$.
The mapping $\pi _{n+1}$ is a surjective homomorphism
of graded rings, i.e
\begin{eqnarray}
\pi _{n+1} : \Lambda _{n+1} \to \Lambda _{n},\;\;\forall r\in\mathbb{N},\;\;
{\pi _{n+1}}^r :={\pi _{n+1}} _{\mid _{{\tiny\Lambda _{n+1}^r}}} :
\Lambda _{n+1}^r \to \Lambda _{n}^r.
\nonumber
\end{eqnarray}
The mapping $\pi _{n+1}^r$ is surjective $\forall r \geq 0$
and an isomorphism if and only if $r\leq n$.
\end{lemma}
This Lemma can be generalized as follows.
\begin{cor}\label{util} Let $n$ be a nonnegative integer.
For any $p \in \mathbb{N},\; p\neq 0$,
the mapping $\Pi_{n+p}: \Lambda_{n+p} \to \Lambda_{n}$,
defined by setting ${\bf x}_{n+1}=0,\;\;{\bf x}_{n+2}=0,\;\dots, {\bf x}_{n+p}=0$,
is a surjective homomorphism of graded rings. Furthermore, the restriction
\begin{eqnarray}
\Pi_{n+p}\mid_{\Lambda_{n+p}^r}:=\Pi_{n+p}^r : \Lambda_{n+p}^r \to \Lambda_{n}^r
\end{eqnarray}
is surjective for all $r\geq 0$, and an isomorphism if and only if $r\leq n$.
\end{cor}
In the following, the notation $A \equiv B$ means that the set $A$ is in
bijection with $B$. Note that, here, since the group homomorphism
(linearity) is insured, group bijection means group isomorphism. So,
in the following, we will use one or other terminology to refer to the same property.

{\bf Proof of Corollary \ref{util}.}
We proceed by induction on $p$. The order $p=1$ corresponds to
Lemma \ref{sim}, i.e $\;\;\Pi_{n+1}^r \equiv \pi_{n+1}^r$ and
$\;\;\Lambda_{n+1}^r \equiv \Lambda_{n}^r$. The surjectivity is then immediate
$\forall p \in \mathbb{N},\; p\neq 0$, as $r\geq 0$. For the {\em one to one} property,
suppose the statement holds for the order $p-1$. For $p$, setting $n+p-1= n'$ and
$n+p= n'+1$
and using Lemma \ref {sim}, $\Lambda_{n'+1=n+p}^{r}\equiv
\Lambda_{n'=n+p-1}^{r}\;\;\Leftrightarrow
\;\;r \leq n+p-1$. Proceeding step by step, we get
$ \Lambda_{n+1}^{r}\equiv \Lambda_{n}^{r}\;\;\Leftrightarrow \;\;r
\leq n;\;\; \Lambda_{n+2}^{r}\equiv \Lambda_{n+1}^{r}\;\;\Leftrightarrow\;\;
r\leq n+1;\dots;
\Lambda_{n+p}^{r}\equiv \Lambda_{n+p-1}^{r}\;\;\Leftrightarrow \;\;r \leq n+p-1.$
Therefore, $\Lambda_{n+p}^{r}\equiv\Lambda_{n}^{r}$ requires $r \leq \;min\;
(n+p-1\;,\dots,\;n+1,\;n)=n$. This ends the proof of the corollary.\hfill$\square$

\begin{exmp} Given $n=2$ and $r= 2$ so that the set of
indeterminates is $\left\{ {\bf x}_1, {\bf x}_2 \right\}$, then,
the following polynomials
$f_i$ are symmetric and of degree $2$, i.e belong to $\Lambda_{2}^{2}$:
\begin{eqnarray}
f_1= {\bf x}_1 {\bf x}_2, \;\; f_2= {\bf x}_1^2 + {\bf x}_2^2, \;\;
 \forall p,q \in \mathbb{Z},\;\;
f= pf_1 + qf_2 \in \Lambda_{2}^{2}.
\nonumber
\end{eqnarray}
Adding a new indeterminate ${\bf x}_3$,
we have the corresponding elements
of $\Lambda_{3}^{2}$
\begin{eqnarray}
&&f'_1 = {\bf x}_1 {\bf x}_2 + {\bf x}_1 {\bf x}_3
+ {\bf x}_2 {\bf x}_3, \; f'_2= {\bf x}_1^2 + {\bf x}_2^2 + {\bf x}_3^2,\cr
&& \forall p,q \in\mathbb{Z},\;  f'= pf'_1 +qf'_2 \in \Lambda_{2}^{2}.
\nonumber
\end{eqnarray}
with $\pi_{3} f'_{i} =f_{i}$, for $i=1,2.$
\end{exmp}
From Lemma \ref{sim}, the following statement holds.
\begin{cor} The sequence of groups
$ 0 \stackrel{i}{\longrightarrow} \Lambda_{n+1}^r
\stackrel{\pi_{n+1}}{\longrightarrow}
\Lambda_{n}^r \stackrel{p}{\longrightarrow} 0$, where $i$ is the canonical injection,
and $p$ the projection onto $\left\{ 0 \right\}$, is exact if and only
if $r\leq n$.
\end{cor}
This corollary may be of great importance for $r>n$
in the Homology Theory involving the groups of
symmetric polynomials \cite{bh}.
\begin{defn}\label{limit}
Let $r$ be a nonnegative integer.
The projective (or inverse) limit $\;\; \Lambda ^{r}
=\lim_{\stackrel{\longleftarrow}{n}}\Lambda _{n}^r
$ is the additive group of sequences of
homogeneous symmetric polynomials of degree $r$ such that
$f^{r}= (f_{1}^{r},f_{2}^{r},\dots , f_{n}^{r}, \dots)$ with
\begin{eqnarray}\label{lim}
\forall n\in \mathbb{N}\not\,\;\left\{0 \right\},\;\;
f_{n}^{r} \in \Lambda_{n}^r\;\;\mbox{and}\;\;
 \pi _{n+1}(f_{n+1}^{r}) = f_{n}^{r}.
\end{eqnarray}
The elements of $\Lambda ^{r}$ are called projective limits and
 $\Lambda ^{r}$ is called the homogeneous group of degree $r$ of projective
limits. Besides, let $\; \Lambda =
\bigoplus _{r\geq 0} \Lambda ^{r} $ be the graded ring defined by  the direct
sum of the homogeneous groups $\Lambda ^{r}$. An element $f$ of $\Lambda$ is a sum
of projective limits, namely $f = \sum _{r \geq 0} f^r $ such that,
 for any degree $r$, $f^r$ belongs to the homogeneous group $\Lambda ^{r}$.
 An element of $\Lambda$ is called a symmetric function.
\end{defn}
It can be shown the following statement \cite{IGM}.
\begin{prop}\label{granpi}
 With the above notation,
 there is a surjective homomorphism of graded rings
 ${\bf {\Pi}}_{n}: \Lambda \to \Lambda_{n}$ defined
 by setting ${\bf x}_{p}=0, \forall p\geq n+1$.
\end{prop}
\begin{exmp}
 Given two nonnegative integers $r$ and $n$, the partial sum $
f_{n}^{r}= \sum_{i=1}^{n} {\bf x}_{i}^r $ defines the sequence
$\left( f_{n}^{r}\right)_{n\in \mathbb{N}}$ as a projective limit of $\Lambda ^{r}$.
This symmetric
function is of degree $r$ and defined by $f^{r}= \sum_{n=1}^{\infty} {\bf x}_{n}^r$.
\end{exmp}
\begin{remark}
Often in the literature, there is no distinction between the projective
limit which is a sequence and the limit of the corresponding partial
sum which is a function. In any case, given $f_n$, the expression
$ \lim_{\stackrel{\longleftarrow}{n}} f_n^r = f^r$ contains all information
generated by the equation {\rm (\ref{lim})}.
\end{remark}
More rigorously, we consider also the following definition
of the inverse limit \cite{L}.

Let $I$ be a set of indices. Suppose a given relation of partial ordering in $I$.
We say that $I$ is directed if given $i,j \in I$, there is $k \in I$ such that
$i\leq k$ and $j\leq k$. Assume that $I$ is directed.
Let now consider $\mathfrak{A}$ a category, and $\left\{ A_i \right\}$
a family of objects in
$\mathfrak{A}$. For each pair $i,j$ such that $i\leq j$, let us consider a given
 morphism:
\[ f_{(j,i)}: A_j\to A_i \]
such that, whenever $i\leq k\leq j$,
one gets
\[
f_{(j,k)}\circ f_{(k,i)}= f_{(j,i)} \;\;\mbox{and}\;\; f_{(i,i)}=id,
\]
where $id$ is the identity mapping of $A_i$.
Such a family is called a directed family of morphisms.
An inverse limit for the family $(f_{(j,i)})$ is a universal object
of the following category $\mathfrak{C}$. $Ob(\mathfrak{C})$ consists of pairs
$(A, (f_i) )$ where $A\in Ob(\mathfrak{A})$ and $(f_i)$ is a family of
morphisms $f_i: A\to A_i$, $i\in I$, such that, for all $i \leq j$, the following
diagram is commutative:
\begin{center}
\unitlength=1mm
\begin{picture}(20,30)(30,35)
\put(37.85,60){$A$}
\put(25,52){$f_j$}
\put(37,60){\vector(-1,-1){15}}
\put(17,40){$A_j$}
\put(42,60){\vector(1,-1){15}}
\put(51,52){$f_i$}
\put(23,40){\vector(1,0){33}}
\put(57,40){$A_i$}
\put(37.25,37){$f_{(j,i)}$}
\end{picture}
\end{center}
Given two nonnegative integers $n_1 \leq n_2$,
let $({\Pi}_{(n_2,n_1)})$  be
the family of homomorphisms of graded rings from $\Lambda_{n_2}$ to $\Lambda_{n_1}$
such that
\[
 {\Pi}_{(n_2,n_1)}:= \Pi_{n_1+(n_2-n_1)},
\]
where $\Pi_{n}$ is defined by Corollary \ref{util}.
We can easily check that $({\Pi}_{(n_2,n_1)})$ is a
directed family of ring homomorphisms in the category of graded rings.
The family $({\bf \Pi}_n)$ of Proposition \ref{granpi}
defines the following commutative diagram:
$\forall m,n \in \mathbb{N},\;\; m\geq n,$
\begin{center}
\unitlength=1mm
\begin{picture}(20,30)(30,35)
\put(37.85,60){$\Lambda$}
\put(22,52){${\bf{\Pi}}_m$}
\put(37,60){\vector(-1,-1){15}}
\put(17,40){$\Lambda_{m}$}
\put(42,60){\vector(1,-1){15}}
\put(51,52){${\bf {\Pi}}_n$}
\put(24,40){\vector(1,0){33}}
\put(58,40){$ \Lambda_{n}$}
\put(37.25,36){$\Pi_{(m,n)}$}
\end{picture}
\end{center}
 $(\Lambda, {\bf\Pi}_{n})$, considered as a universal object,
is the inverse limit of the directed family $({\Pi}_{(n_2,n_1)})$.
So, we agree with the property that
{\em the inverse limit defined by the family
of directed homomorphisms $(\Pi_{(n_2,n_1)})$ is equal to the inverse limit
defined by the family of projection $(\pi_{n})$.}
Here and thereafter, defining inverse limit by the family
 $\{\Lambda_{n}, \Pi_{(n_1,n_2)} \}$
or by the family $\{\Lambda_{n}, \pi_{n}\}$ is {\em equivalent}.

\section{Multi-indicial symmetric functions}
In this section, we define the symmetric function of infinite number
of entries that we call multi-indicial symmetric function.

Given $m,n,k \in \mathbb{N}$, let us consider the following set of independent
indeterminates
\[
\left\{ {{\bf a}}_{ \left\{ 1\leq p \leq m,
\left[ \mu \right]_{1\leq i\leq k}\right\} } \right\}.\]
See Table \ref{tab1}.
The notation $\left[ \mu \right]_{i}$ means any multi-index of the form
$\mu_1\mu_2\dots\mu_i$,
with $1\leq \mu_i \leq n$, $1\leq i\leq k$. For instance,
${\bf a}_{m\left[ \mu \right]_{p}}$
 denotes in general
${\bf a}_{m\mu_1\mu_2\dots\mu_p}$, for any $1\leq\mu_i\leq n$.
The sets of indeterminates may be organized in the
following manner:
\begin{eqnarray}
&&D^{0}_{m}=\{{\bf a}_m\},\; D^{0}_{(m)} = \bigcup_{1\leq l \leq m}D^{0}_{l},\;
 D^{k}_{m} = \bigcup _{1 \leq \mu_1,\mu_2,\dots,\mu_k \leq n} {\bf a}
_{m\mu_1\mu_2\dots\mu_k},\;\\
&& D^{k}_{(m)} =\bigcup _{l=1}^{m} D^{k}_{l},
 D^{(k)}_{m}=    \bigcup _{l=0}^{k} D^{l}_{m},\;
\mathcal D_{(m)} = \bigcup _{k = 0}^{\infty} D^{k}_{(m)},\;\; \mathcal D^{(k)} =
\bigcup _{m = 1}^{\infty} D^{(k)}_{m}, \\
&&\mathcal D = \mathcal D_{(m)}\cup\mathcal D^{(k)}.
\end{eqnarray}
The following statement holds by a simple combinatorics.
\begin{prop} Let $m$, $n$ and $k$ be three nonnegative integers.
Then,
\[
\; \mid D^{(k)}_{m+1} \mid = q_k,\;\;\mid D^{k+1}_{(m)} \mid=n^{k+1}m,\;\;
 \mid {\mathcal D}^{(k+1)}_{(m+1)} \mid = q_k + n^{k+1}(m+1),
 \]
where  $\ {\mathcal D}^{(k+1)}_{(m+1)}=D^{(k)}_{m+1}\cup D^{k+1}_{(m)}\cup
\left\{{\bf a}_{(m+1)\mu_1\mu_2\dots\mu_{k}\mu_{k+1}}\right\}$,
$q_k=\frac{n^{k+1}-1}{n-1}$ if $n\neq 1$ and $q_k= k+1$ if $n=1$.
\end{prop}
\begin{defn}
Given $m$, $n$, $k$, three nonnegative integers such that $\;m,\; n\geq 1$, then
the polynomial ring
\begin{eqnarray}
&&\mathbb{Z} ({\bf a}_1, {\bf a}_2, \dots, {\bf a}_m,{\bf a}_{1\mu},{\bf a}_{2\mu}, \dots
, {\bf a}_{m\mu},
 {\bf a}_{1\mu_{1}\mu_{2}},{\bf a}_{2\mu_{1}\mu_{2}},\dots,
{\bf a}_{m\mu_{1}\mu_{2}},\dots,\\
&& \;\;\;\,{\bf a}_{1\mu_{1}\mu_{2}\dots\mu_{k}},
{\bf a}_{2\mu_{1}\mu_{2}\dots\mu_{k}},\dots,
{\bf a}_{m\mu_{1}\mu_{2}\dots\mu_{k}}
), \nonumber
\end{eqnarray}
where $\mu,\;\mu_{1},\;\mu_{2},\;\;\dots,\; \mu_{k-1}$ and $\;\mu_{k}$
take all values
in  $\left\{\, 1,2,\dots, n\,\right\}$,
is denoted by $\mathbb{Z} \left[{\bf a}_m{\bf a}_{m{\left[ \mu \right] }_k}\right]$.
The number of indeterminates is $mq_k$, where
$
q_k = \frac{n^{k+1}- 1}{n-1}
\;\;{\mbox{if}}\;\; n\neq 1; \;\;
q_k =  k+1, \;{\mbox{if}}\;\; n = 1.
$
The symmetric group $S_{mq_k} $
defines the graded ring of symmetric polynomials of
$\mathbb{Z} [ {\bf a}_m {\bf a}_{m\left[\mu\right]_k}]$:
$\Lambda_{m,k} = \mathbb{Z}[ {\bf a}_m
{\bf a}_{m{\left[\mu\right]}_k}] ^{S_{mq_k}}. $
\end{defn}
\begin{lemma}\label{base} Let $r$, $m$ and $k$ be two nonnegative integers.
There is a group isomorphism $\;\; \Lambda_{m,k}^r \equiv\Lambda_{m q_k}^r \;$
leading to a graded ring isomorphism $\; \Lambda_{m,k} \equiv \Lambda_{m q_k} \;$
\end{lemma}
{\bf Proof.} The set of indeterminates
\begin{eqnarray}
\left\{\left\{ {\bf a}_p \right\}, \left\{ {\bf a}_{p\mu}\right\}_{p; \mu=1,...,n},
\dots,
\left\{ {\bf a}_{p\mu_1\mu_2\dots\mu_{k}}\right\}_{p; \mu_1,\mu_2,\dots,\mu_{k}= \;1,
...,n}
\right\}_{p=1,...,m}
\nonumber
\end{eqnarray}
can be viewed as the set of
indeterminates $\; \left\{ {\bf x}_{1}, {\bf x}_{2}, \dots, {\bf x}_{m q_k}\right\}. \;$
The independence of indeterminates requires the
${\bf a}_{m\mu_{1}\mu_{2}\dots\mu_{p}}$ to correspond to a unique
${\bf x}_{i}$. Since the two sets possess the same cardinal, a well defined
bijection can be built from one onto the other.
\hfill$\square$

By convention,
$\Lambda_{m,0}=\Lambda_{m},\;\;\Lambda_{0,k}=\Lambda_{q_{k-1}},\;\;
\Lambda_{0,0}=\mathbb{Z}.$

\begin{defn}
Let $m$ and $k$ be nonnegative integers, $m\geq 1$.
\begin{enumerate}
\item[(i)]
The graded ring homomorphism
$\;\;h_{m+1,k}: \Lambda_{m+1, k}\to\Lambda_{m, k},\;$
such
that $\forall r\in \mathbb{N},\;\,$
$h_{m+1,k}\mid_{\Lambda_{m+1,k}^{r}}:=h_{m+1,k}^r:\Lambda_{m+1, k}^{r}\to\Lambda_{m,
k}^{r}\;$ and defined by setting
\[
{\bf a}_{m+1}=0,\; {\bf a}_{(m+1)\mu}=0,\; \dots,
{\bf a}_{(m+1)\mu_{1}\mu_{2}...\mu_{k}}=0,
\;\;\forall\; 1\leq\mu_i\leq n,
\]
is called the $(m+1,k)$ horizontal projection or simply the h-projection
when no confusion occurs;
the restriction $h_{m+1,k}^r$ is called the $(m+1,k)$ horizontal projection
of degree $r$.
\item[(ii)]
The graded ring homomorphism $\; v_{m,k+1}:\Lambda_{m, k+1}\to\Lambda_{m,k},\;$
such that
$\;\forall r\in \mathbb{N},\;\;$
$\;v_{m,k+1}\mid_{\Lambda_{m,k+1}^{r}}:=v_{m,k+1}^r: \Lambda_{m, k+1}^{r}
\to \Lambda_{m, k}^{r}\;$
and defined by setting
\begin{equation}
{\bf a}_{1\mu_{1}\mu_{2}...\mu_{k+1}}=0, \;\;
{\bf a}_{2\mu_{1}\mu_{2}...\mu_{k+1}}=0,\;\;
{\bf a}_{m\mu_{1}\mu_{2}...\mu_{k+1}}=0,
\;\;\forall\; 1\leq\mu_i\leq n,
\nonumber
\end{equation}
is called $(m,k+1)$ vertical
projection or simply the v-projection when no confusion occurs; the
restriction $v_{m,k+1}^r$ is called the $(m,k+1)$ vertical projection of
degree $r$.
\item[(iii)]
The graded ring homomorphism $\;\pi_{m+1,k+1}:\Lambda_{m+1,
k+1}\to\Lambda_{m, k},\;$ such that $\forall r\in\mathbb{N},\;\;$
$\pi_{m+1,k+1}\mid_{\Lambda_{m+1,k+1}^{r}}:
=\pi_{m+1,k+1}^r:\Lambda_{m+1,k+1}^{r}\to\Lambda_{m, k}^{r}\;$
and defined by setting
\begin{eqnarray}
&&{\bf a}_{m+1}=0,\;
{\bf a}_{(m+1)\mu}=0,\; \dots, \,{\bf a}_{(m+1)\mu_{1}\mu_{2}...\mu_{k}}=0,\;
{\bf a}_{1\mu_{1}\mu_{2}...\mu_{k+1}}=0,\cr
&&{\bf a}_{2\mu_{1}\mu_{2}...\mu_{k+1}}=0,\;
{\bf a}_{m\mu_{1}\mu_{2}...\mu_{k+1}}=0\;\;\mbox{and}\cr
&&{\bf a}_{(m+1)\mu_1\mu_2\dots\mu_{k+1}}=0,\;
\forall\; 1\leq\mu_i\leq n,
\nonumber
\end{eqnarray} is called the
$(m+1,k+1)$ projection. The restriction $\pi_{m+1,k+1}^r$ is called the
$(m+1,k+1)$ projection of degree $r$.
\end{enumerate} \end{defn}

\begin{lemma}\label{crucial} Given $m$ a nonnegative integer, $m\geq 1$, the
h-projection $ h_{m+1,1}^r:\Lambda_{m+1,1}^{r} \to \Lambda_{m,1}^{r}$
is surjective for all $r \geq 0$ and bijective if and only if $r \leq (n+1)m$.
\end{lemma}
{\bf Proof.} The results follow from Lemma \ref{base} and
Corollary \ref{util}. The surjectivity is immediate. For the proof of the
bijectivity, we obtain using Lemma \ref{base}
$\Lambda_{m+1,1}^{r} \equiv \Lambda_{(m+1)q_1}^r$ and
$\Lambda_{m,1}^{r} \equiv \Lambda_{mq_1}^r$.
From Corollary \ref{util}, the r.h.s expressions are bijective
if and only if $0\leq r\leq mq_1$. \hfill$\square$
\begin{prop}\label{surjecfin}
\begin{enumerate}
\item[(i)] $\;\;\forall k\geq 0,\;$ the h-projection
$h_{m+1,k}^r: \Lambda_{m+1,k}^{r} \to \Lambda_{m,k}^{r}$ is surjective for all $r
\geq 0$ and bijective if and only if $r \leq mq_k$.
\item[(ii)] $\;\;\forall m\geq 1,\;$ the v-projection $ v_{m,k+1}^r: \Lambda_{m,k+1}^{r} \to
\Lambda_{m,k}^{r} $ is surjective for all $r \geq 0$ and bijective if and only
if $r \leq mq_k$.
\item[(iii)] $\;\;\forall k\geq 0,\;\;\forall m\geq 1$, the projection
$\pi_{m+1,k+1}^r: \Lambda_{m+1,k+1}^{r} \to \Lambda_{m,k}^{r} $ is surjective for
all $r \geq 0$ and bijective if and only if $r \leq mq_k$.
\end{enumerate}
\end{prop}
{\bf Proof.} The proofs of the surjections are immediate
by the use of Lemma \ref{base}. So, let us pay attention to the proofs of
the bijections. One can show $(i)$ by induction on $k$. Consider Lemma
\ref{crucial} as the order $k=1$. The following step is similar to the
proof of Corollary \ref {util}, taking the $min$ on different values of
$r \leq \;min\; \left\{ (m+1)q_k, m q_k \right\} = mq_k$. The steps $(ii)$ and
$(iii)$ can be shown by the same way. Indeed, we can easily give the
prescribed equivalent of Lemma \ref{crucial} for $k$ and for both $m$
and $k$.
\hfill$\square$

\begin{defn} Given $r$ a nonnegative integer, then
\begin{enumerate}
\item[(i)] We call horizontal (resp. vertical) sequence (m,k) of degree $r$
the inverse system denoted by $(\Lambda_{m,k}^{r}, h^{r}_{m,k})_{m\in\mathbb{N}}$,
(resp. $(\Lambda_{m,k}^{r}, v^{r}_{m,k})_{k\in\mathbb{N}}$);
\item[(ii)] We call sequence (m,k) of degree $r$ the inverse system 
denoted by $(\Lambda_{m,k}^{r},
\pi^{r}_{m,k})_{m,k\in \mathbb{N}}.$
\end{enumerate}
\end{defn}
The next proposition can be deduced from Proposition \ref{surjecfin}.
\begin{prop}\label{cocomut} With
the above notation, the following diagram in which
all mappings are surjective for all $r\in \mathbb{N}$, and bijective if $r\leq
mq_k$, is commutative
\begin{center}
\unitlength=1mm
\begin{picture}(65,55)(0,-5)
\put(-5,40){$\Lambda_{m+1,k+1}^{r}$}
\put(12,40){\vector(1,0){42}}
\put(29,42){$h^{r}_{m+1,k+1}$}
\put(12,35){\vector(3,-2){42}}
\put(33,22){$\pi_{m+1,k+1}$}
\put(56,40){$\Lambda^{r}_{m,k+1}$}
\put(0,35){\vector(0,-1){28}}
\put(-15,22){$v^{r}_{m+1,k+1}$}
\put(12,3){\vector(1,0){42}}
\put(30,0){$h^{r}_{m+1,k}$}
\put(-5,3){$\Lambda_{m+1,k}^{r}$}
\put(62,35){\vector(0,-1){28}}
\put(63,22){$v_{m,k+1}^{r}$}
\put(57,3){$ \Lambda_{m,k}^{r}$}
\end{picture}
\end{center}
and leads to the corresponding commutative diagram with respect to the
graded ring structure.
\end{prop}
The previous development leads to the following consequence.
Given a nonnegative integer $r$, taking the projective limit with respect to the
horizontal (resp. vertical)  sequence $(m,k)$ of degree $r$,
we obtain $\tilde\Lambda^{r}_{.,k} = \lim_{\stackrel{\longleftarrow}{m}}
\Lambda^{r}_{m,k}$
(resp. $\tilde\Lambda^{r}_{m,.} = \lim_{\stackrel{\longleftarrow}{k}}
\Lambda^{r}_{m,k}$ ) that we call the horizontal (resp. vertical)
projective limit of degree $r$ of the sequence
$(\Lambda_{m,k}^{r})_{m\in \mathbb{N}}$ (resp. $(\Lambda_{m,k}^{r})_{k\in \mathbb{N}}$).
$\tilde\Lambda^{r}_{.,k}$ (resp. $\tilde\Lambda^{r}_{m,.}$) is the additive
group of horizontal (resp. vertical) symmetric functions of degree $r$.
Furthermore, given $m\in\mathbb{N}$ (resp. $k\in\mathbb{N}$), for each $k$ (resp. $m$),
there is a surjective homomorphism
\begin{eqnarray}
H_{m,k}^{r}: \Lambda_{.,k}^{r}\to \Lambda_{m,k}^{r} \;\;
(\;\;\mbox{resp.}\;\;V_{m,k}^{r}: \Lambda_{m,.}^{r}\to \Lambda_{m,k}^{r}\;\;)
\nonumber
\end{eqnarray}
defined by ${\bf a}_{p>m}=0$ and $\forall q\in \mathbb{N},$
${\bf a}_{p>m; \left[\mu \right]_{0\leq q\leq k}}=0$
(resp. ${\bf a}_{1\leq p\leq m; \left[\mu\right]_{ q\geq (k+1)}}$ $=0$)
 which is bijective
iff. $r\leq mq_k$.
\begin{remark} The elements of
$\;\;\tilde\Lambda^{r}_{m,.}$ and $\;\;\tilde\Lambda^{r}_{.,k}$ are sequences of
symmetric functions of a given degree $r$. The groups
 $\;\;\tilde \Lambda^{r}_{m,.}$ and $\tilde\Lambda^{r}_{.,k}$ are not isomorphic.
Indeed, one way to easily realize this is to notice that
the ring $\mathbb{Z}\left[ {\bf a}_{m}{\bf a}_{m \left[  \mu \right]_k}\right]$
has not the same dependence with respect to the indeterminates
${\bf a}_{m}$ and ${\bf a}_{m \left[\mu \right]_k}$.
Implicitly, ${\bf a}_{m \left[\mu \right]_k}$ depends on $n$, while,
in an obvious manner,
${\bf a}_{m}$ does not. Thus the elements of $\Lambda^{r}_{m,.}$, at the limit
$k\to \infty$, do not involve the integer parameter $n$ at the opposite
of those of $\tilde\Lambda^{r}_{.,k}$ as $m\to \infty$.
This construction of the ring $\mathbb{Z} \left[ {\bf a}_{m}{\bf a}_{m \left[ \mu \right]_k}\right]$ is different
from the construction of a polynomial ring in the indeterminates {\rm \cite{bh}}
$\left\{ {\bf a}_{m,n,k}\right\}_{0\leq m \leq M;\; 0\leq n \leq N; 0\leq k \leq K}$,
given $M,N, K \in \mathbb{N}$, which consists in assigning the free entries of a
 $3$-tensor,
for instance.
The independence between the indeterminates, in this case,
should correspond to the isomorphism of
sets of sequences of symmetric functions in the remaining indices.
We have
 $\Lambda_{m,n,(.)}^{r}\equiv\Lambda_{m,(.),k}^{r}\equiv\Lambda_{(.),n,k}^{r}$,
 where the point means that
 the corresponding index tends to infinity.
\end{remark}
Summing over the degrees, one obtains the graded rings
$\tilde\Lambda_{.,k} = \bigoplus_{r \geq 0}\tilde\Lambda^{r}_{.,k}$,
$\;\;\tilde\Lambda_{m,.}=\bigoplus_{r \geq 0}\tilde\Lambda^{r}_{m,.}$
 of horizontal sequences of symmetric functions
and vertical sequences of symmetric functions, respectively.
\begin{defn} Let $r$ be a nonnegative integer.
Two symmetric
functions $P^{r}$ and $Q^{r}$ of degree $r$ are said equal
if and only if
$\;\forall n \in \mathbb{N},\;\; P^{r}_{n}=Q^{r}_{n}$.
\end{defn}
Let $P^{r}_{m,.} \in \;\tilde\Lambda^{r}_{m,.}.$
 For any $k\in
\mathbb{N},\;\;\;P^{r}_{m,.}=(P^{r}_{m,0},P^{r}_{m,1},\dots,P^{r}_{m,k},\dots)$
such that, for any $k$, $v^{r}_{m,k+1}P^{r}_{m,k+1}=P^{r}_{m,k}$.
We also obtain, for any $m \in \mathbb{N}$, $h^{r}_{m+1,k}P^{r}_{m+1,k}=P^{r}_{m,k}$.
Hence, the mapping $h_{m}^{r}:\tilde\Lambda^{r}_{m,.} \to \tilde\Lambda^{r}_{m-1,.}$
defined by
\begin{eqnarray}
h_{m}^{r}(P^{r}_{m,.})=(h_{m,0}^{r}P^{r}_{m,0},h_{m,1}^{r}P^{r}_{m,1},
\dots ,h_{m,k}^{r}P^{r}_{m,k},\dots),
\nonumber\end{eqnarray}
allows to get
$h_{m}^{r}(P^{r}_{m,.})=P^{r}_{m-1,.}$. This shows that $h_{m}^{r}$ is a
well defined projection and defines the projective limit of degree $r$ of
the vertical sequence
$\left(\Lambda_{m,.}^{r}\right)_{m\in \mathbb{N}}$ by ${\mathcal P}^{r}
=\lim_{\stackrel{\longleftarrow}{m}}P^{r}_{m,.}$.
We call this inverse limit the h(v)-limit of degree $r$. Besides, defining the
mapping $v^{r}_{k}:\tilde\Lambda^{r}_{.,k}\to \tilde\Lambda^{r}_{.,k}$, by
\begin{eqnarray}
v^{r}_{k}(P^{r}_{.,k})=(v^{r}_{0,k}P^{r}_{0,k},v^{r}_{1,k}P^{r}_{1,k},
\dots, v^{r}_{m,k}P^{r}_{m,k}, \dots),
\nonumber\end{eqnarray}
we get $v_{k}^{r}(P^{r}_{.,k})=P^{r}_{.,k-1}$ which shows that
$v_{k}^{r}$ is a well
defined projection which defines the projective limit of degree $r$ of the
horizontal sequence $(\Lambda_{.,k}^{r})_{k\in \mathbb{N}}$
by ${\mathcal P'}^{r}=\lim_{\stackrel{\longleftarrow}{k}}P^{r}_{.,k}$.
We call this inverse limit the v(h)-limit of degree $r$.
${\mathcal P}^{r}$ and ${\mathcal P'}^{r}$ are not {\em a priori} the same quantity.
But, they are actually isomorphic. Indeed from \cite{L},
the following holds.
\begin{theo}
\label{doubl}
Let M and K be two directed sets, $\left( A_{m,k}\right)_{m\in M; k\in K}$
be a family of Abelian groups equipped with homomorphisms
labeled by $M\times K$, and defining an inverse limit.
Assigning the obvious ordering to the product $M\times K$,
i.e
$(m,k)\leq (m',k')\; \Leftrightarrow\; m\leq m'$ and $k\leq k',$
the following inverse limits exist and are isomorphic in a natural way:
\begin{eqnarray}
\label{dirdedir}
\lim_{\stackrel{\longleftarrow}{m}}\lim_{\stackrel{\longleftarrow}{k}}A_{m,k}=
\lim_{\stackrel{\longleftarrow}{k}}\lim_{\stackrel{\longleftarrow}{m}}A_{m,k}.
\end{eqnarray}
\end{theo}
The inverse systems $((\Lambda_{m,k}^{r}, v_{m,k}^{r}), h_{m}^{r})$ and
$((\Lambda_{m,k}^{r}, h_{m,k}^{r}), v_{k}^{r})$,
giving rise to the v(h)-limit and the h(v)-limit, respectively, are equivalent.
We establish this equivalence in the following.
\begin{prop}\label{dede}
Let $r,k,m\;$ ($m\geq 1$) be three nonnegative integers.
For all $m_1,m_2,$ $k_1,k_2\; \in \mathbb{N}$, such that
$1\leq m_1\leq m_2$ and $0\leq k_1 \leq k_2$,
the mappings
\begin{eqnarray}
\phi_{(m_2,m_1),k}^{r}: \Lambda_ {m_2,k}^{r}\to \Lambda_ {m_1,k}^{r} \;\;
\mbox{and}\;\;
\psi_{m, (k_2,k_1)}:\Lambda_ {m,k_2}^{r}\to \Lambda_ {m,k_1}^{r},
\nonumber\end{eqnarray}
defined by $\phi_{(m_1,m_1),k}^{r}\equiv\mathbb{I}$,
$\psi_{m, (k_1,k_1)}\equiv\mathbb{I}$,
\begin{eqnarray}
\phi_{(m_2,m_1),k}^{r}\equiv
h_{(m_1+1),k}^{r}\circ h_{(m_1+2),k}^{r}\circ\dots\circ h_{m_2,k}^{r}
\nonumber\end{eqnarray}
and
\begin{eqnarray}
\psi_{m, (k_2,k_1)}^{r}\equiv
v_{m,(k_1+1)}^{r}\circ v_{m,(k_1+2)}^{r}\circ\dots\circ v_{m,k_2}^{r}
\nonumber\end{eqnarray}
are well defined surjective group homomorphisms.
Furthermore, given $k$ (resp. $m$),
 $(\phi_{(m_2,m_1),k})$ (resp. $(\psi_{m,(k_2,k_1)})$) defines a directed family
 of homomorphisms of graded rings.
\end{prop}
{\bf Proof.}
The surjectivity is given by induction from the definition of
$h_{m,k}^{r}$ and $v_{m,k}^{r}$. Moreover, one can easily check, that given $k$,
for any $m_1\leq p\leq m_2,$
$\phi_{(m_2, p),k}^{r}\circ \phi_{(p, m_1),k}^{r} =\phi_{(m_2, m_1),k}^{r}.$
Given m, the similar property also holds for $\psi_{m,(k_2,k_1)}^{r}$.
\hfill$\square$
\begin{lemma}\label{phipsi}
With the above notation,
\begin{eqnarray}
&&\phi_{(m,m-1),k}^{r}\equiv h_{m,k}^{r},\;\;
\psi_{m, (k,k-1)}^{r}\equiv v_{m,k}^{r},\;\;
\phi_{(m,m-1),k}^{r}\circ\psi_{m, (k,k-1)}^{r}=\pi_{m,k}^{r},\nonumber\\
\label{commuta}
&& \phi_{(m_2,m_1),k_1}^{r}\circ\psi_{m_2,(k_2,k_1)}^{r}
 \equiv \psi_{m_1,(k_2,k_1)}^{r}\circ\phi_{(m_2,m_1),k_2}^{r},\\
&&\forall 0\leq q\leq k_2,\;\;\phi_{(m_2,m_1),q}^{r}\circ\phi_{(m_2,m_1),k_2}^{r}
 =\phi_{(m_2,m_1),k_2}^{r},\cr
&&\forall 1\leq p\leq m_2,\;\;\psi_{m_2,(k_2,k_1)}^{r}\circ\psi_{p,(k_2,k_1)}^{r}
=\psi_{m_2,(k_2,k_1)}^{r}.
\end{eqnarray}
Furthermore, one has, $\forall m_1\leq p \leq m_2$ and
$\forall k_1\leq q \leq k_2$,
\begin{eqnarray}
\phi_{(m_2,p), k}^{r}\circ\phi_{(p,m_1), k}^{r}= \phi_{(m_2,m_1), k}^{r},
 \;\;\mbox{and}\;\;
\psi_{m, (k_2,q)}^{r}\circ\psi_{m, (q,k_1)}^{r}= \psi_{m, (k_2,k_1)}^{r}
\nonumber\end{eqnarray}
\end{lemma}
{\bf Proof.}
 This is immediate from Propositions \ref{cocomut} and \ref{dede}.
\hfill$\square$
\begin{remark}
(\ref{commuta}) can be viewed as the data of a commutative diagram.
\end{remark}
Given nonnegative integers $k,m\geq 1$, the directed families
$(\phi_{(m_2,m_1),k}^{r})$ and $(\psi_{m, (k_2,k_1)}^{r})$
define a directed family of homomorphisms
in both the indices $m$ and $k$ as follows.
$\forall\, m_1, m_2, k_1, k_2$ $\in \mathbb{N}$,
such that $1\leq m_1\leq m_2$ and $0\leq k_1 \leq k_2$,
let
\begin{eqnarray}
\Phi^{r}_{(m_2,m_1),(k_2,k_1)}: \Lambda_{m_2,k_2}^{r} \to \Lambda_{m_1,k_1}^{r}
\nonumber\end{eqnarray}
be the mapping defined by
\begin{eqnarray}
&&\Phi^{r}_{(m_1,m_1),(k_2,k_1)}\equiv\psi_{m_1, (k_2,k_1)}^{r},\;\;
\Phi^{r}_{(m_2,m_1),(k_1,k_1)}\equiv\phi_{(m_2,m_1),k_1}^{r},\\
&&\Phi^{r}_{(m_2,m_1),(k_2,k_1)}\equiv \phi_{(m_2,m_1),k_1}^{r}
\circ
\psi_{m_2, (k_2,k_1)}^{r}=
\psi_{m_1, (k_2,k_1)}^{r}\circ\phi_{(m_2,m_1),k_1}^{r}.
\end{eqnarray}
We deduce, from Lemma \ref{phipsi},
with $m_1\leq p\leq m_2$ and $k_1\leq q\leq k_2$,
\begin{eqnarray}
\label{granphi}
&&\Phi^{r}_{(m_2,p),(k_2,k_1)}\circ\Phi^{r}_{(p,m_1),(k_2,k_1)}
\equiv\Phi^{r}_{(m_2,m_1),(k_2,k_1)},\\
&&\Phi^{r}_{(m_2,m_1),(k_2,q)}\Phi^{r}_{(m_2,m_1),(q,k_1)}
\equiv\Phi^{r}_{(m_2,m_1),(k_2,k_1)}.
\nonumber
\end{eqnarray}
Applying Theorem \ref{doubl} with $M= \mathbb{N} \not\;\{0\}$ and
$K=\mathbb{N}$ which are obviously directed sets,
\begin{eqnarray}
\left(\Phi^{r}_{(m_2,m_1),(k_2,k_1)}
\right)_{m_{1,2}\in M; k_{1,2}\in K}
\nonumber\end{eqnarray}
is a directed family of homomorphisms labeled by
$M\times K$ which allows  to write,
 by analogy with (\ref{dirdedir})
 \begin{eqnarray}
\lim_{\stackrel{\longleftarrow}{m}}
\lim_{\stackrel{\longleftarrow}{k}}\Lambda_{m,k}^{r}
=\lim_{\stackrel{\longleftarrow}{k}}
\lim_{\stackrel{\longleftarrow}{n}}\Lambda_{m,k}^{r}.
\nonumber
\end{eqnarray}
The following statement is valid.
\begin{theo}
\begin{enumerate}
\item[(i)]Let $r,m_1,m_2,k_1,k_2$ be nonnegative integers such that
 $1 \leq m_1\leq m_2$
and $0\leq k_1\leq k_2$.
The mappings $\phi_{(m_2,m_1)}: \Lambda_{m_2,.} \to \Lambda_{m_1,.}$
and $\psi_{(k_2,k_1)}: \Lambda_{.,k_{2}}\to\Lambda_{.,k_{1}}$ defined by
\begin{eqnarray}
&&\phi_{(m_2,m_1)}^r=
h_{m_{1}+1}^{r}\circ h_{m_{1}+2}^{r}\circ \dots\circ h_{m_{2}}^{r},\cr
&&\psi_{(k_2,k_1)}^{r}=v_{k_{1}}^{r}\circ v_{k_{1}+1}^{r}\circ
 \dots\circ v_{k_{2}}^{r}
\nonumber
\end{eqnarray}
are surjective homomorphisms of graded rings, define directed families
$(\phi_{(m_2,m_1)})$ and $(\psi_{(k_2,k_1)})$ of homomorphisms of graded rings.
Furthermore, the
inverse limits induced by these families are equal, i.e
\begin{eqnarray}
\lim_{\stackrel{\longleftarrow}{m}}\Lambda_{m,.}=\Lambda
=\lim_{\stackrel{\longleftarrow}{k}}\Lambda_{.,k}.
\nonumber
\end{eqnarray}
\item[(ii)] Given a nonnegative integer $m$ (resp. $k$), $m>0$,
the mapping $H_{m}: \Lambda \to \Lambda_{m,.}$
(resp. $V_{k}: \Lambda \to\Lambda_{.,k}$) defined by setting ${\bf a}_{p> m}=0$
and $\;\forall q\in \mathbb{N},\;\; {\bf a}_{p>m;\left[\mu\right]_{q}}=0$
(resp. ${\bf a}_{p>0}=0$ and ${\bf a}_{p>0;\left[\mu \right]_{q> k}}=0$)
is a surjective homomorphism of graded rings.
\end{enumerate}
\end{theo}
{\bf Proof.}
 We prove that the two inverse limits coincide.
The statement, mainly obtained by the definition of any universal object
of a category,
holds in general by Theorem \ref{doubl}. Let us illustrate, here, this statement
by a particular case. We consider that $k$ and $m$ are two nonnegative
integers with $m>0$. Moreover, homomorphism means surjective homomorphism
of graded rings. Given two nonnegative integers $m>0$ and $k$,
there are four directed families of homomorphisms
\begin{eqnarray}
 (\phi_{(m_1,m_2),k}),\;\; (\psi_{m,(k_1,k_2)}),\;\; (\phi_{(m_1,m_2)})\;\;
 \mbox{and}\;\;
(\psi_{(k_1,k_2)}),
\nonumber\end{eqnarray}
generating four kinds of categories $\mathcal{C}_{k}$, $\mathcal{C}_{m}$,
$\mathfrak{C}_{1}$ and $\mathfrak{C}_{2}$ whose the sets of objects are given by
\begin{eqnarray}
&& Ob(\mathcal{C}_{k})= \{ (\mathcal{R},\; (\mathcal{H}_{m,k})_{m} )\},\;\;
\mathcal{H}_{m,k}: \mathcal{R} \to \Lambda_{m,k},\label{ob1}
\\
&& Ob(\mathcal{C}_{m})= \{(\mathcal{R},\;  (\mathcal{V}_{m,k})_{k} )\},\;\; \mathcal{V}_{m,k}:
\mathcal{R} \to \Lambda_{m,k},\\
&& Ob(\mathfrak{C}_{1})=\{ (\mathcal{R},\; (\mathcal{H}_{m})_{m} )\},\;\; \mathcal{H}_{m}:
 \mathcal{R}\to \Lambda_{m,.},\\
&& Ob(\mathfrak{C}_{2})=\{ (\mathcal{R},\; (\mathcal{V}_{k})_{k} )\},\;\; \mathcal{V}_{k}:
 \mathcal{R}\to \Lambda_{k,.}, \label{ob4}
\end{eqnarray}
respectively, where $\mathcal{R}$ is any graded ring.
The categories of (\ref{ob1})-(\ref{ob4}) generate,
up to a unique isomorphism, universal objects given by
\begin{eqnarray}
 (\Lambda_{.,k}, (H_{m,k})_{m}),\;\;(\Lambda_{m,.}, (V_{m,k})_{k}),\;\;
(\tilde{\Lambda}_{1}, (H_{m})_{m})\;\;\;\mbox{and}\;\;\;
(\tilde{\Lambda}_{2}, (V_{k})_{k}),
\nonumber\end{eqnarray}
respectively. Let us consider the following diagram
\begin{center}
\unitlength=1mm
\begin{picture}(60,50)(-20,0)
\put(-25,40){$\tilde{\Lambda}_{2}^{r}$ }
\put(-7,40){ \vector(-1,0){10} }
\put(-12,49){ $f$ }
\put(-17,46){\vector(1,0){10} }
\put(-12,43){ $g$}
\put(-5,40){ $\tilde{\Lambda}_{1}^{r}$ }
\put(5,40){\vector(1,0){30} }
\put(19,42){$H_{m}^{r}$}
\put(5,35){\vector(3,-2){30}}
\put(22,25){$A_{m,k}$}
\put(38,40){$\Lambda^{r}_{m,.}$}
\put(-22,35){\vector(0,-1){20}}
\put(-27,25){$V^{r}_{k}$}
\put(-15,10){\vector(1,0){50}}
\put(10,5){$H^{r}_{m,k}$}
\put(-25,10){$\Lambda_{.,k}^{r}$}
\put(40,35){\vector(0,-1){20}}
\put(41,25){$V_{m,k}^{r}$}
\put(37,10){$ \Lambda_{m,k}^{r}$}
\end{picture}
\end{center}

Given $k$, for all $m$, it comes $A_{m,k}:=V_{m,k}\circ H_{m}:
 \tilde{\Lambda}_1 \to \Lambda_{m,k}$.
Moreover, $(\tilde{\Lambda}_1, (A_{m,k})_{m})\in  Ob(\mathcal{C}_{k})$
 and there is a unique homomorphism
$\varphi_{k}:\tilde{\Lambda}_1 \to \Lambda_{.,k}$ such that:
\begin{eqnarray}
V_{m,k}\circ H_{m}= H_{m,k} \circ\varphi_{k}.
\nonumber\end{eqnarray}
It follows that $(\tilde{\Lambda}_1, (\varphi_{k})_{k} )\in Ob(\mathfrak{C}_{2})$ and
thus, there exists a unique homomorphism
$g: \tilde{\Lambda}_1 \to \tilde{\Lambda}_2 $ such that
$\varphi_{k} = V_{k}\circ g$. Hence, we get
\begin{eqnarray}\label{funct1}
 V_{m,k}\circ H_{m}= H_{m,k} \circ V_{k}\circ g.
\end{eqnarray}
Moreover, in the same manner, given $m$ and the homomorphism
$H_{m,k}\circ V_{k}: \tilde{\Lambda}_2 \to \Lambda_{m,k}$, for all $k$,
we have, through the universal object property of $(\Lambda_{m,.}, (V_{m,k})_{k})$,
the unique homomorphism
$\varpi_{m}: \tilde{\Lambda}_2 \to \Lambda_{m,.}$ such that
\begin{eqnarray}
H_{m,k}\circ V_{k}= V_{m,k}\circ \varpi_{m}.
\nonumber\end{eqnarray}
$\varpi_{m}$ induces, by the universal object property of
$(\tilde{\Lambda}_{1}, (H_{m}))$,
the factorization $\varpi_{m,2}=H_{m}\circ f$, where
$f:\tilde{\Lambda}_{2}\to\tilde{\Lambda}_{1}$
is uniquely defined. Consider
\begin{eqnarray}\label{funct2}
H_{m,k}\circ V_{k}= V_{m,k}\circ H_{m}\circ f.
\end{eqnarray}
From (\ref{funct1}) and (\ref{funct2}), we deduce  $f\circ g= \mathbb{I}$.
Conversely, we can also show that $g\circ f=\mathbb{I}$.
\hfill$\square$
\begin{theo}
The diagram defined by the Table \ref{tab2}
is commutative in the sense that any of its squares is commutative.
\end{theo}
{\bf Proof.}
Given four nonnegative integers $m,k,p,q$, with $m>0$,
any of the {\em internal} diagrams, i.e any diagram of the form
\begin{center}
\unitlength=1mm
\begin{picture}(65,50)(0,0)
\put(-5,40){$ \Lambda_{m+p,k+q}^{r} $}
\put(12,40){\vector(1,0){40}}
\put(24,42){$\phi^{r}_{(m+p, m),k+q}$}
\put(12,35){\vector(3,-2){40}}
\put(26,25){$\Phi_{(m+p,m),(k+q,k)}$}
\put(55,40){$\Lambda^{r}_{m,k+q}$}
\put(0,35){\vector(0,-1){27}}
\put(-20,22){$\psi^{r}_{m+p,(k+q,k)}$}
\put(12,4){\vector(1,0){40}}
\put(25,1){$\phi^{r}_{(m+p,m),k}$}
\put(-5,4){$\Lambda_{m+p,k}^{r}$}
\put(59,35){\vector(0,-1){27}}
\put(60,22){$\psi_{m,(k+q,k)}^{r}$}
\put(56,4){$ \Lambda_{m,k}^{r}$}
\end{picture}
\end{center}
is commutative from Proposition \ref{phipsi} and the properties (\ref{granphi}).
Let us pay attention to the diagrams involving the inverse limits.
There are three kinds of such diagrams.
\begin{enumerate}
\item[(i)]
The first involves two inverse limits in $m$:
\begin{center}
\unitlength=1mm
\begin{picture}(65,50)(0,0)
\put(5,40){$\Lambda_{.,k_2}^{r}$}
\put(15,40){\vector(1,0){35}}
\put(25,42){$\psi^{r}_{(k_2,k_1)}$}
\put(52,40){$\Lambda^{r}_{.,k_1}$}
\put(7,35){\vector(0,-1){15}}
\put(-6,30){$H^{r}_{m,k_2}$}
\put(17,15){\vector(1,0){33}}
\put(25, 10){$\psi^{r}_{m,(k_2,k_1)}$}
\put(5,15){$\Lambda_{m,k_2}^{r}$}
\put(53,35){\vector(0,-1){15}}
\put(54,30){$H^{r}_{m,k_1}$}
\put(52,15){$ \Lambda_{m,k_1}^{r}$}
\end{picture}
\end{center}
Such a diagram is commutative for $(\Lambda_{.,k_1}, (H_{m,k_1}))$ is a
universal object.
$\psi_{m,(k_2,k_1)}\circ H_{m,k_2}: \Lambda_{.,k_2}\to \Lambda_{m,k_1}$
can be factorized
by the unique ring homomorphism
 $\Lambda_{.,k_2}\stackrel{\psi_{(k_2,k_1)}}{\longrightarrow}\Lambda_{.,k_1}$
as
$\psi_{m,(k_2,k_1)}\circ H_{m,k_2}=H_{m,k_1}\circ\psi_{(k_2,k_1)}.$
\item[(ii)]
The second involves two inverse limits in $k$:
\begin{center}
\unitlength=1mm
\begin{picture}(65,50)(0,0)
\put(5,40){$\Lambda_{m_2,.}^{r}$}
\put(15,40){\vector(1,0){35}}
\put(25,42){$\phi^{r}_{(m_2,m_1)}$}
\put(51,40){$\Lambda^{r}_{m_1,.}$}
\put(7,35){\vector(0,-1){15}}
\put(-6,27){$V^{r}_{m_2,k}$}
\put(16,15){\vector(1,0){34}}
\put(25,10){$\phi^{r}_{(m_2,m_1),k}$}
\put(5,15){$\Lambda_{m_2,k}^{r}$}
\put(53,35){\vector(0,-1){15}}
\put(54,27){$V_{m_1,k}^{r}$}
\put(52,15){$ \Lambda_{m_1,k}^{r}$}
\end{picture}
\end{center}
The diagram is also commutative by the usual definition of
universal object $(\Lambda_{m_{1},.}, (V_{m_{1},k}))$ by analogy with the
 proof of the case $(i)$.
\item[(iii)]
The third involves the inverse limit $\Lambda$:
\begin{center}
\unitlength=1mm
\begin{picture}(65,50)(0,0)
\put(5,40){$\Lambda^{r}$}
\put(12,40){\vector(1,0){37}}
\put(24,42){$H^{r}_{m}$}
\put(51,40){$\Lambda^{r}_{m,.}$}
\put(7,35){\vector(0,-1){15}}
\put(0,27){$V^{r}_{k}$}
\put(12,15){\vector(1,0){37}}
\put(25,10){$H^{r}_{m,k}$}
\put(5,15){$\Lambda_{.,k}^{r}$}
\put(53,35){\vector(0,-1){15}}
\put(58,27){$V_{m,k}^{r}$}
\put(52,15){$ \Lambda_{m,k}^{r}$}
\end{picture}
\end{center}
\end{enumerate}
The commutativity results from the same argument.
\hfill$\square$

The commutativity of any of the diagrams defined by Table \ref{tab2}
represents the inverse limit defined by the three directed families
$(\phi_{(m_1,m_2)})$, $(\psi_{(k_1,k_2)})$ and  $(\Phi_{(m_2,m_1)(k_2,k_1)})$.
Thus,
\begin{eqnarray}
\left\{ \left\{ \Lambda_{m,k}, v_{m,k} \right\} , h_m \right\}\Leftrightarrow
\left\{ \left\{ \Lambda_{m,k}, h_{m,k} \right\} , v_k \right\}\Leftrightarrow
\left\{ \Lambda_{m,k}, \Phi_{(m_2,m_1),(k_2,k_1)} \right\}
\nonumber\end{eqnarray}
that leads to
\begin{eqnarray}
 \tilde\Lambda =\lim_{\stackrel{\longleftarrow}{m}}
\lim_{\stackrel{\longleftarrow}{k}}\Lambda_{m,k}^{r}
=\lim_{\stackrel{\longleftarrow}{m}}\tilde\Lambda_{m,.} =\lim_{\stackrel{\longleftarrow}{k}}\tilde\Lambda_{.,k}.
 \end{eqnarray}
Finally, the set of usual symmetric functions is recovered, i.e.
 $ \tilde\Lambda \equiv \Lambda.$
\section{Multi-partition and multi-indicial Schur functions}
In this section, we deal with the definition of multi-partition and study
the corresponding interesting family of symmetric functions
known as the Schur functions \cite{IGM}.
\begin{defn}\label{partition}
A partition $\lambda$ is a finite or infinite sequence of integers
$(\lambda_{1},\lambda_{2},$ $\dots,\lambda_{i}, \dots)$, with
$\lambda_{1}\geq \lambda_{2}\geq\dots\geq 0$ and $\mid\lambda\mid:=\sum_{i}\lambda_{i}<\infty,$
so that,
from a certain point onwards (if $\lambda$ is infinite), all the $\lambda_{i}$ are
$0$. The non zero $\lambda_i$ are called the parts of $\lambda$. The number of
parts is the length $l(\lambda)$ of $\lambda$.
\end{defn}
\begin{remark}
Two partitions $\lambda_1$ and $\lambda_2$ which differ only by a sequence of $0$ at
the end are equal. For instance, $(1,2)$ and $(1,2,0,0,\dots)$ are regarded
as the same partition.
\end{remark}
\begin{defn} Given two nonnegative
integers $m$ and $k$, we call a $\left[ m,k \right]$-partition (or a
multi-partition when no confusion occurs) the ordered sequence
\begin{eqnarray}
\lambda_{\left[
m,k \right]}= (\lambda^{\left[ m,0 \right]}, \lambda^{\left[ m,1 \right]},
\dots,\lambda^{\left[ m,k \right]})
\nonumber\end{eqnarray}
defined by a set of $k+1$ partitions such that:
\begin{eqnarray}\label{multiparty}
&&\lambda^{\left[ m,0 \right]}=
\lambda^{\left[ m,0 \right]}_{1\leq i} =(\lambda^{\left[ m,0 \right]}_1,\lambda^{\left[ m,0
\right]}_2,\dots,\lambda^{\left[ m,0 \right]}_p,\dots),\;\;
\lambda^{\left[ m,1 \right]}=(\lambda^{\left[ m,1 \right]}_{p\mu})_{1\leq p,\mu},\cr
&&\lambda^{\left[ m,k \right]}= ({\lambda^{\left[ m,k
\right]}}_{p\mu_1\mu_2\dots\mu_k})_{ 1\leq p,\mu_1,\mu_2,\dots,\mu_k},
\end{eqnarray}
with $1\leq \mu_{l}\leq n,\;\;for\;\;1\leq l\leq k$, so that the following property
is satisfied:
\begin{eqnarray}\label{ineq}
&&\lambda^{\left[ m,0 \right]}_{1}\geq\lambda^{\left[ m,0 \right] }_{2}\geq\dots\geq
\lambda^{\left[ m,0\right] }_{p}\geq\dots\geq\lambda^{\left[ m,0 \right]}_{m}\geq
\dots\geq\lambda^{\left[ m,1\right]}_{11}\geq\cr
&&\lambda^{\left[ m,1\right]}_{12}\geq\dots\geq\lambda^{\left[ m,1 \right] }_{1\mu}\geq\dots
\geq\lambda^{\left[ m,1 \right]}_{21}\geq\dots\geq\lambda^{\left[ m,1 \right] }_{2\mu}\geq\dots
\geq\lambda^{\left[ m,1\right]}_{2n}\geq\dots\cr
&&\geq\lambda^{\left[ m,1\right]}_{p\mu}\geq\dots
\lambda^{\left[ m,1\right]}_{m1}\geq\lambda^{\left[ m,1\right]}_{m2}\geq\dots
\geq\lambda^{\left[ m,1\right]}_{mn}\geq\dots \geq\lambda^{\left[ m,k\right]}_{m\mu_1\dots \mu_k}\cr
&&\geq\lambda^{\left[ m,k\right]}_{m\mu_1\dots (\mu_k+1)}\geq\dots\geq 0.
\end{eqnarray}
$\lambda^{\left[ m,p \right]}$, for any $0 \leq p \leq k$,
is called a sub-partition
of $\lambda_{\left[ m,k \right]}$.
Furthermore, we identify
\begin{eqnarray}
\mid\lambda_{\left[ m,k \right]}\mid=\sum_{0 \leq
p\leq k} \mid\lambda^{\left[ m,p\right]}\mid.
\nonumber
\end{eqnarray}
 The length of the
$\left[ m,k \right]$-partition is defined by the
sum of the lengths of its
sub-partitions, namely
\begin{eqnarray}
l(\lambda_{\left[ m,k \right]})=\sum_{0 \leq p\leq k}
l(\lambda^{\left[ m,p \right]}).
\nonumber\end{eqnarray}
\end{defn}
One can easily see that the so defined $\left[ m,k \right]$-partition is 'exhaustive'
relatively to the number of
indeterminates, i.e it assigns an exponent to each of them. Furthermore, a
$\left[ m,0 \right]$-partition is, by convention, a partition in the sense of Definition \ref{partition}.
Each of the $\lambda^{\left[ m, k\right]}$, taken separately, defines a partition
such that the ordered sequence (\ref{multiparty}) which defines $\lambda_{\left[ m,k \right]}$ remains a
partition. One can define the {\em monomial symmetric function} in
the $mq_k$ indeterminates by the sum of all distinct monomials that can be obtained
from
\begin{eqnarray}
\quad
{\bf a}^{\lambda_{\left[ m,k \right]}}= {\bf a}_{1}^{\lambda^{\left[ m,0
\right]}_1}\dots {\bf a}_{m}^{\lambda^{\left[ m,0 \right]}_m}
\prod_{m;\mu}{\bf a}_{m\mu}^{\lambda^{\left[
m,1 \right]}_{m\mu}} \dots \prod_{m;\mu_1,\mu_2,\dots,\mu_k}
{\bf a}_{m\mu_{1}\mu_{2}\dots\mu_{k}}^{{\lambda^{\left[ m,k
\right]}}_{m\mu_1\mu_2\dots\mu_k}},
\label{monosym}
\end{eqnarray}
by permutation of the ${\bf a}$'s.
In particular, for any  $i\in \left[ 0,k \right]$,
\begin{eqnarray}
{\lambda_{\left[ m,k \right]}}=(
\lambda^{\left[ m,0 \right]}=(0),\dots,\lambda^{\left[ m,i\right]}
=(\underbrace{1,1,\dots,1}_{r-times},0,0, \dots), \dots, \lambda^{\left[ m,k \right]}=(0)).
\nonumber\end{eqnarray}
One readily recovers the definition of
classical symmetric monomial $e_{1}=m_{(1^{r})}$ \cite{IGM}. It is then immediate
that $\mathbb{Z}$-basis of $\tilde \Lambda_{.,k}$ and $\tilde\Lambda_{m,.}$ can be obtained
as a function of the monomial symmetric functions corresponding to
(\ref{monosym}), when the $\left[ m,k \right]$-partition runs through all
multi-partitions.
Let us come back to the usual theory.
Let $n$ be a nonnegative integer. In the following,
$\delta$ is the partition defined by $(n-1,n-2,\dots,1,0)$.
The following statement holds \cite{IGM}.
\begin{prop}
Given a nonnegative integer $n$,
for each partition $\alpha = (\alpha_{1},$ $\alpha_{2},\dots,\alpha_{n})$,
of nonnegative integers such that $\alpha_{1}>\alpha_{2}>\dots >\alpha_{n}\geq 0,$
the homogeneous polynomial defined by
\begin{eqnarray}\label{aalph}
a_{\alpha}=\rm{det}({\bf x}_{i}^{\alpha_{j}})_{1\leq i,j \leq n}
\end{eqnarray}
is divisible by the Vandermonde determinant $a_{\delta}$ in
$\mathbb{Z}\left[{\bf x}_1, {\bf x}_2,\dots, {\bf x}_n \right]$.
\end{prop}
The partition $\alpha$ can be chosen as $\alpha_{i}=\lambda_{i}+ (n-i)$, for $1\leq
i\leq n$, so that $\alpha=\lambda + \delta$, where $\lambda$ is a partition of length
at most $n$. The quotient $s_{\lambda}({\bf x}_1, \dots,{\bf x}_n)= a_{\lambda + \delta}\not \;
a_{\delta}$ is a symmetric polynomial, homogeneous of degree $\mid \lambda \mid$.
Passing to $n+1$ variables, we have
\begin{eqnarray}
s_{\lambda}({\bf x}_1, \dots,{\bf x}_n,
{\bf x}_{n+1})\mid _{{\bf x}_{n+1}=0}= s_{\lambda}({\bf x}_1, \dots,{\bf x}_n,0)=s_{\lambda}({\bf x}_1,
\dots,{\bf x}_n).
\nonumber\end{eqnarray}
 The uniquely defined quotient $s_{\lambda} \in \Lambda$, that reduces to
$s_{\lambda}({\bf x}_1, \dots,{\bf x}_n)$ when ${\bf x}_{p\geq n+1}=0$, for any $n\geq
l(\lambda)$, is the {\em Schur function} corresponding to the partition $\lambda$.
Let us consider the $mq_k$ indeterminates with $n\geq 1$
 (see Table \ref{tab1}). We define the
 $\left[ m,k \right]$-partition
\begin{eqnarray} &&\delta_{\left[ m,k \right]}=
 \left(\delta^{\left[ m,0 \right]}, \delta^{\left[ m,1\right]},\dots ,\delta^{\left[ m,k \right]}\right),
 \;\;\mbox{by}\;\;
 \delta^{\left[ m,0 \right]}= (\delta^{\left[ m,0 \right]}_p=mq_k-p)_{1\leq p\leq m},\cr
 &&\delta^{\left[ m,1\right]}=
 (\delta^{\left[ m,1 \right]}_{p\mu}= mnq_{k-1}-(p-1)n-\mu)_{(1\leq p\leq
m);(1\leq\mu\leq n)},
\nonumber\end{eqnarray}
and, for any $0\leq d\leq k$, $1\leq p\leq m$ and
$1\leq\mu_i\leq n$,
\begin{eqnarray}
\delta^{\left[ m,d\right]}_{p\mu_1\mu_2\dots\mu_d}&=&\delta^{\left[
m,d-1\right]}_{mn\dots n}-(p-1)n^{d}-\sum_{1\leq l\leq
d-1}(\mu_l-1)n^{d-l}-\mu_d \cr
&=&n^{d}(mq_{k-d}-(p-1))-\sum_{1\leq l\leq
d-1}(\mu_l-1)n^{d-l}-\mu_d,
\nonumber\end{eqnarray}
where the index $mn\dots n$ contains
$(d-1)$ times the index $n$.
Explicitly, it can be written
\begin{eqnarray}
\delta^{\left[ m,0 \right]}&=&(
\delta^{\left[ m,0 \right]}_1=mq_k-1,\delta^{\left[ m,0\right]}_2=mq_k-2,\dots,
\delta^{\left[ m,0\right]}_p=mq_k-p,\dots,\cr
&&\delta^{\left[ m,0 \right]}_m=mnq_{k-1}),
\nonumber\end{eqnarray}
where the identity $q_k-1= nq_{k-1}$ has been used.
\begin{eqnarray}
\delta^{\left[ m,1 \right]}&=&( \delta^{\left[
m,1 \right]}_{11} =mnq_{k-1}-1,\;\;\delta^{\left[ m,1 \right]}_{12}=mnq_{k-1}-2,\dots,\cr
&&\delta^{\left[ m,1 \right]}_{1\mu}=mnq_{k-1}-\mu,\;\;
\delta^{\left[ m,1 \right]}_{1n}=n(mq_{k-1}-1), \cr
&&\delta^{\left[ m,1 \right]}_{21}=n(mq_{k-1}-1)-1,\dots,
\delta^{\left[ m,1\right]}_{2n}=n(mq_{k-1}-2),\cr
&&\dots,\delta^{\left[ m,1 \right]}_{p\mu}=
mnq_{k-1}-(p-1)n-\mu,\dots, \cr
&& \delta^{\left[ m,1 \right]}_{mn}= mn(q_{k-1}-1)=mn^2q_{k-2} ),\cr
\delta^{\left[ m,2 \right]}&=&(
\delta^{\left[ m,2 \right]}_{1,1,1}=mn^2q_{k-2}-1, \dots,\cr
&&\delta^{\left[ m,2 \right]}_{p\mu_1\mu_2}=mn^2q_{k-2}-(p-1)n^2-(\mu_1 - 1)n-\mu_2, \dots, \cr
&&\delta^{\left[ m,2 \right]}_{mnn}= mn^3q_{k-3} ),\dots \nonumber\end{eqnarray}

Finally, $\delta^{\left[ m,k \right]}_{111\dots
1}=mn^k-1$ and
 $\delta^{\left[ m,k \right]}_{mnn\dots n}=0$. Hence, $\delta_{\left[ m,k \right]}$
realizes a partition such that
\begin{eqnarray}
a_{\delta_{\left[ m,k\right]}}= {\rm det}\left(
{\bf a}_{p\mu_1\mu_2\dots\mu_d}^{\delta^{\left[ m,t\right] }_{q\nu_1\nu_2\dots\nu_t}} \right)_{(1\leq
p,q\leq m);(0\leq t,d\leq k);(1\leq \mu_i, \nu_j\leq n)}
\nonumber\end{eqnarray}
corresponds to the Vandermonde determinant of the matrix $\mathcal A$,
see Table \ref{tab3}.
The following statement holds.
\begin{prop} Let
$\lambda_{\left[ m,k \right]}$ be a multi-partition of length $l(\lambda_{
\left[ m,k \right]})\geq mq_{k}$ such that the inequalities  {\rm(\ref{ineq})} are
strict. There exists a $\left[ m,k\right]$-partition $\ell_{\left[ m,k \right]}$ of
length at most $mq_k$ such that
\begin{eqnarray}\label{multcond}
&&\forall 1\leq p\leq m,\; \forall 0\leq d\leq k, \forall 1\leq\mu_1,\mu_2,\dots,\mu_d
\leq n,\,\,\cr
&&\lambda^{\left[ m,d \right]}_{p\mu_1\mu_2\dots\mu_d} = \ell^{\left[ m,d
\right]}_{p\mu_1\mu_2\dots\mu_d} +
\left\{ mq_k -(p-1)n^d - \left(\sum_{l=1}^{d}(\mu_{l}-1)n^{d-l}+1\right) \right\}
\end{eqnarray}
\end{prop}
{\bf Proof.}
One has to consider the previous construction
of $\delta_{\left[ m,k\right]}$ whose length is $mq_k-1$.
Each of the sub-partitions of $\delta_{\left[ m,k\right]}$,
namely the $\delta^{\left[ m,d\right]}$, has a length less than or equal to the
length of the sub-partition $\lambda^{\left[ m,d \right]}$ of
$\lambda_{\left[ m,k \right]}$.
 Indeed, if the inequalities (\ref{ineq}) are strict,
 for any $0\leq d\leq k$, $1\leq p\leq m$,
 \begin{eqnarray}
 \ell^{\left[ m,d
\right]}_{p\mu_1\mu_2\dots\mu_d}= \lambda^{\left[ m,d\right]}_{p\mu_1\mu_2\dots\mu_d}-
\delta^{\left[m,d\right]}_{p\mu_1\mu_2\dots\mu_d}\geq 0,
\nonumber\end{eqnarray}
that allows to define the
$\left[ m,k \right]$-partition $\ell_{\left[ m,k \right]}$ whose length
is at most $mq_k$.
\hfill$\square$

We write
$\lambda_{\left[ m,k\right]}= \delta_{\left[ m,k\right]}+
\ell_{\left[ m,k \right]}$. In an obvious manner,
$a_{\lambda_{\left[ m,k \right]}}$ is divisible by the Vandermonde determinant
$a_{\delta_{\left[ m,k \right]}}$. The quotient
\begin{eqnarray}
S_{\ell_{\left[ m,k \right]}}=
a_{\lambda_{\left[ m,k \right]}}\not\;\; a_{\delta_{\left[ m,k \right]}}
\nonumber\end{eqnarray}
is of course a symmetric
polynomial, homogeneous of degree $\mid \ell_{\left[ m,k \right]}\mid$.
For the sake of simplicity, let us set $\mid \ell \mid:=\mid \ell_{\left[ m,k \right]}\mid$.
Passing to $m+1$ (resp. $k+1$), we set
$a_{\lambda_{\left[ m+1,k \right]}}\not\;\; a_{\delta_{\left[ m+1,k \right]}}= S_{\ell_{\left[ m+1,k\right]}}$,
(resp. $\;a_{\lambda_{\left[ m,k+1 \right]}}\not\;
a_{\delta_{\left[ m,k+1 \right]}}= S_{\ell_{\left[ m,k+1 \right]}}$).
Under the horizontal (resp. vertical)
projection $h_{m+1,k}^{\mid \ell\mid}$ (resp. $v_{m,k+1}^{\mid \ell\mid}$),
 we have $S_{\ell_{\left[ m+1,k \right]}}\to
S_{\ell_{\left[ m,k \right]}}$ (resp. $S_{\ell_{\left[ m,k+1 \right]}}\to
S_{\ell_{\left[ m,k \right]}}$).
This horizontal (resp. vertical) sequence defines an unique horizontal (resp. vertical)
inverse limit $S_{\ell_{\left[ .,k \right]}}\in \Lambda_{.,k}$
(resp. $S_{\ell_{\left[ m,. \right]}} \in \Lambda_{m,.}$) which reduces to $S_{\ell_{\left[ m,k \right]}}$
setting $\forall p\geq 1, \forall k\geq 0,
{\bf a}_{m+p}=0,\dots, {\bf a}_{(m+p)\left[\mu\right]_{k}}=0$
(resp. $\forall m\geq 0, \forall p\geq
1, {{\bf a}_{m}}_{\left[ \mu\right]_{k+p}}=0$), for any $mq_k \geq
l(\ell_{\left[ m,k \right]})$. We call $S_{\ell_{\left[ .,k \right]}}$
(resp. $S_{\ell_{\left[ m,. \right]}}$)
the horizontal (resp. vertical) Schur function corresponding to $\ell$.
Taking the inverse limit with respect to the other index,
 one recovers the usual Schur function corresponding
to $\ell$.

\section{Concluding remarks}
In this paper, we have extended
 the symmetric functions to the multi-indicial symmetric functions.
 The multi-indicial symmetric functions can be viewed as
 the elements of the universal objects in the sense of the inverse limits
 of the categories ${\mathcal{C}}_{m;k}$ (\ref{ob1})-(\ref{ob4}).
  By the construction of
 the corresponding partition, called multi-partition,
 we have defined the multi-indicial Schur functions.
 Further properties of multi-indicial symmetric functions as well as
 the relation to the symmetric functions $P_{\lambda}(q,t)$ \cite{IGM}
 will be discussed in the forthcoming paper \cite{bh}.

\section*{Acknowledgments}
M.N.H. thanks the National Institute for Theoretical Physics (NITheP) and its Director Prof. Frederik G. Scholtz for hospitality during a pleasant stay in Stellenbosch where
a part of this work has been finalized.
This work was supported under a grant of
the  National Research Foundation of South Africa
and by the ICTP through the OEA-ICMPA-Prj-15.
The ICMPA is in partnership with the Daniel Iagolnitzer Foundation (DIF), France.

\section*{Appendix: Examples of multi-indicial Schur polynomials with $n=2$}
$van_{mk}$ denotes the Vandermonde determinant
and $spoly_{mk}$ the Schur polynomial associated with $\ell_{\left[ m,k \right]}$.

{\bf Example 1: $m=1,\;\;$ $k=1$}
\begin{eqnarray}
 {\ell}_{\left[ 1,1 \right]} &=&\left[ \begin {array}{cc} {\ell}^{[1,0]}= \left[
\begin {array}{c} 3\end {array} \right], &{\ell}^{[1,1]}= \left[
\begin {array}{cc} 2,&1\end {array} \right] \end {array} \right];  \cr
{ van_{11}}&&= \left( Y_{{1,1}}-Y_{{1,2}} \right)  \left( X_{{1}}-Y_{{1,
2}} \right)  \left( X_{{1}}-Y_{{1,1}} \right); \cr
{spoly_{11}}&=&Y_{{1,1}}X_{{1}}Y_{{1,2}} \left( {X_{{1}}}^{2}+Y_{{1,2}}
X_{{1}}+Y_{{1,1}}X_{{1}}+{Y_{{1,2}}}^{2}+{Y_{{1,1}}}^{2}+Y_{{1,2}}Y_{{
1,1}} \right).
\nonumber
\end{eqnarray}

{\bf Example 2: $m=2,\;\;$$k=1$}
\begin{eqnarray}
{\ell}_{\left[ 2,1\right]}&=& \left[ \begin {array}{cc} {\ell}^{[2,0]}= \left[
\begin {array}{cc} 3,&2\end {array} \right], &{\ell}^{[2,1]}= \left[
\begin {array}{cccc} 2,&1,&1,&1\end {array} \right] \end {array} \right]; \cr
van_{21}&=& \left( -Y_{{2,1}}+Y_{{1,1}} \right)
\left( -Y_{{2,1}}+X_{{2}} \right)
\left( X_{{2}}-Y_{{1,1}} \right)
\left( -Y_{{2,1}}+X_{{1}}  \right) \cr
&& \left( X_{{1}}-Y_{{1,1}} \right)
 \left( X_{{1}}-X_{{2}} \right)
  \times  \left( -Y_{{2,1}}+Y_{{1,2}} \right)
    \left( Y_{{1,1}}-Y_{{1, 2}} \right) \cr
  &&
  \left( -Y_{{1,2}}+X_{{2}} \right)
  \left( X_{{1}}-Y_{{1,2}} \right)
  \left( Y_{{2,1}}-Y_{{2,2}} \right)
   \times \left( -Y_{{2,2}}+Y_{{1,1}} \right) \cr
  &&
    \left( -Y_{{2,2}}+X_{{2}} \right)
   \left( -Y_{{2,2}}+X_{{1}} \right)
 \left( -Y_{{2,2}}+Y_{{1,2}} \right); \cr
spoly_{21}&=& X_{{1}}X_{{2}}Y_{{1,1}}Y_{{1,2}}Y_{{2,1}}Y_{{2,2}}\times  \cr
&&\left( \right. 3\left(\right.
 \,Y_{{2,1}}X_{{2}}Y_{{1,2}}Y_{{2,2}}
+ \,Y_{{2,1}}Y_{{1,2}}X_{{1}}Y_{{1,1}}
+ \,Y_{{2,1}}Y_{{1,2}}X_{{1}}X_{{2}}  \cr
&&
+ \,Y_{{1,2}}X_{{1}}X_{{2}}Y_{{1,1}}
+\,Y_{{1,2}}X_{{2}}Y_{{1,1}}Y_{{2,2}}
+\,X_{{1}}X_{{2}}Y_{{1,1}}Y_{{2,2}} \cr
&&
+\,Y_{{1,2}}X_{{1}}Y_{{1,1}}Y_{{2,2}}
+\,Y_{{2,1}}Y_{{1,2}}X_{{1}}Y_{{2,2}}
+\,Y_{{2,1}}X_{{2}}Y_{{1,2}}Y_{{1,1}}    \cr
&&
+\,Y_{{2,1}}Y_{{1,1}}Y_{{1,2}}Y_{{2,2}}
+\,Y_{{2,1}}X_{{1}}Y_{{1,1}}Y_{{2,2}}
+\,Y_{{2,1}}X_{{2}}Y_{{1,1}}Y_{{2,2}}\cr
&&
+\,Y_{{2,1}}X_{{1}}X_{{2}}Y_{{2,2}}
+\,Y_{{2,1}}X_{{1}}X_{{2}}Y_{{1,1}}
+\,Y_{{1,2}}X_{{1}}X_{{2}}Y_{{2,2}}   \left.\right)
\cr
&&
+  Y_{{1,2}}{X_{{1}}}^{2}X_{{2}}
+ {Y_{{1,2}}}^{2}Y_{{1,1}}Y_{{2,2}}
+ {Y_{{1,2}}}^{2}X_{{2}}Y_{{2,2}}
+ {Y_{{1,2}}}^{2}X_{{2}}Y_{{1,1}} \cr
&&
+ {Y_{{1,2}}}^{2}X_{{1}}Y_{{2,2}}
+ {Y_{{1,2}}}^{2}X_{{1}}Y_{{1,1}}
+ {Y_{{1,2}}}^{2}X_{{1}}X_{{2}}
+ Y_{{1,2}}{Y_{{2,2}}}^{2}Y_{{1,1}} \cr
&&
+ Y_{{1,2}}{Y_{{1,1}}}^{2}Y_{{2,2}}
+ Y_{{1,2}}{Y_{{2,2}}}^{2}X_{{2}}
+ Y_{{1,2}}{Y_{{1,1}}}^{2}X_{{2}}
 +Y_{{1,2}}{X_{{2}}}^{2}Y_{{2,2}} \cr
&&
 +Y_{{1,2}}{X_{{2}}}^{2}Y_{{1,1}}
 +Y_{{1,2}}X_{{1}}{Y_{{2,2}}}^{2}
 +Y_{{1,2}}X_{{1}}{Y_{{1,1}}}^{2}
 +Y_{{1,2}}X_{{1}}{X_{{2}}}^{2} \cr
&&
 +Y_{{1,2}}{X_{{1}}}^{2}Y_{{2,2}}
 +Y_{{1,2}}{X_{{1}}}^{2}Y_{{1,1}}
 +{Y_{{2,1}}}^{2}Y_{{1,1}}Y_{{2,2}}
 +{Y_{{2,1}}}^{2}X_{{1}}X_{{2}}    \cr
&&
 +{Y_{{2,1}}}^{2}X_{{1}}Y_{{2,2}}
 +{Y_{{2,1}}}^{2}X_{{2}}Y_{{2,2}}
+{Y_{{2,1}}}^{2}X_{{1}}Y_{{1,1}}
+{Y_{{2,1}}}^{2}X_{{2}}Y_{{1,1}}  \cr
&&
+{Y_{{2,1}}}^{2}Y_{{1,2}}Y_{{1,1}}
+ {Y_{{2,1}}}^{2}Y_{{1,2}}X_{{2}}
+{Y_{{2,1}}}^{2}Y_{{1,2}}Y_{{2,2}}
 +Y_{{2,1}}{Y_{{1,2}}}^{2}Y_{{2,2}}   \cr
&&
 +Y_{{2,1}}{Y_{{1,2}}}^{2}Y_{{1,1}}
 +Y_{{2,1}}{Y_{{1,2}}}^{2}X_{{2}}
 +Y_{{2,1}}{Y_{{1,1}}}^{2}X_{{2}}
 +Y_{{2,1}}{Y_{{2,2}}}^{2}X_{{2}}  \cr
&&
 +Y_{{2,1}}{Y_{{1,1}}}^{2}Y_{{2,2}}
 +Y_{{2,1}}{Y_{{2,2}}}^{2}Y_{{1,1}}
 +Y_{{2,1}}{X_{{2}}}^{2}Y_{{1,1}}
 +Y_{{2,1}}{X_{{2}}}^{2}Y_{{2,2}}   \cr
&&
 +{Y_{{2,1}}}^{2}Y_{{1,2}}X_{{1}}
 +Y_{{2,1}}Y_{{1,2}}{X_{{1}}}^{2}
 +Y_{{2,1}}{X_{{1}}}^{2}X_{{2}}
 +Y_{{2,1}}{X_{{1}}}^{2}Y_{{1,1}} \cr
&&
 +Y_{{2,1}}{X_{{1}}}^{2}Y_{{2,2}}
 +Y_{{2,1}}X_{{1}}{X_{{2}}}^{2}
 +Y_{{2,1}}X_{{1}}{Y_{{2,2}}}^{2}
 +Y_{{2,1}}X_{{1}}{Y_{{1,1}}}^{2} \cr
&&
 +Y_{{2,1}}Y_{{1,2}}{X_{{2}}}^{2}
 +Y_{{2,1}}{Y_{{1,1}}}^{2}Y_{{1,2}}
 +Y_{{2,1}}{Y_{{2,2}}}^{2}Y_{{1,2}}
 +Y_{{2,1}}{Y_{{1,2}}}^{2}X_{{1}} \cr
&&
 +{X_{{1}}}^{2}X_{{2}}Y_{{1,1}}
 +{X_{{1}}}^{2}X_{{2}}Y_{{2,2}}
 +{X_{{1}}}^{2}Y_{{1,1}}Y_{{2,2}}
 +X_{{1}}{Y_{{1,1}}}^{2}X_{{2}}   \cr
&&
 +X_{{1}}{Y_{{2,2}}}^{2}X_{{2}}
 +X_{{1}}{Y_{{1,1}}}^{2}Y_{{2,2}}
 +X_{{1}}{X_{{2}}}^{2}Y_{{2,2}}
 +{X_{{2}}}^{2}Y_{{1,1}}Y_{{2,2}}\cr
&&
 +X_{{2}}{Y_{{1,1}}}^{2}Y_{{2,2}}
 +X_{{1}}{X_{{2}}}^{2}Y_{{1,1}}
 +X_{{2}}{Y_{{2,2}}}^{2}Y_{{1,1}}
 +X_{{1}}{Y_{{2,2}}}^{2}Y_{{1,1}}
\left.\right).
\nonumber\end{eqnarray}

{\bf Example 3: $m=1,\;\;$$k=2$}
\begin{eqnarray}
{\ell}_{\left[ 1,2\right]}&=& \left[ \begin {array}{ccc} {\ell}^{[1,0]}= \left[
\begin {array}{c} 3\end {array} \right], &{\ell}^{[1,1]}= \left[
\begin {array}{cc} 2,&1\end {array} \right], &{\ell}^{[1,2]}= \left[
\begin {array}{cccc} 1,&1,&1,&1\end {array} \right] \end {array} \right]; \cr
&&\cr
 van_{12}&=& \left( -Z_{{1,2,2}}+Z_{{1,1,2}} \right)
 \left( X_{{1}}-Z_{{1,1,2}} \right)
 \left( X_{{1}}-Z_{{1,2,2}} \right)
\left( Z_{{1,1,1}}-Z_{{1,1,2}} \right) \cr
&\times& \left( Z_{{1,1,1}}-Z_{{1,2,2}} \right)
 \left( -Z_{{1,1,1}}+X_{{1}} \right)
 \left( Y_{{1,2}}-Z_{{1,1,2}} \right)
 \left( Y_{{1,2}}-Z_{{1,2,2}} \right) \cr
&\times&
  \left( X_{{1}}-Y_{{1,2 }} \right)
 \left( Y_{{1,2}}-Z_{{1,1,1}} \right)
\left( Y_{{1,1}}-Z_{{1,1,2}} \right)
\left( Y_{{1,1}}-Z_{{1,2,2}} \right)\cr
&\times&
\left( X_{{1}}-Y_{{1,1}} \right)
\left( Y_{{1,1}}-Z_{{1,1,1}} \right)
\left( Y_{{1,1}}-Y_{{1,2}} \right)
\left( Z_{{1,1,2}}-Z_{{1,2,1}} \right)   \cr
&\times&
 \left( Z_{{1,2,1}}-Z_{{1,2,2}} \right)
\left( -Z_{{1,2,1}}+X_{{1}} \right)
 \left( -Z_{{1,2,1}}+Z_{{1,1,1}} \right)
  \left( -Z_{{1,2,1}} +Y_{{1,2}} \right)\cr
&\times&
   \left( -Z_{{1,2,1}}+Y_{{1,1}} \right);  \cr
spoly_{12}&=&Y_{{1,1}}Y_{{1,2}}Z_{{1,1,1}}Z_{{1,1,2}}Z_{{1,2,1}}Z_{{1,2,2
}}X_{{1}}  \times \cr
&&\left( \right.
 {X_{{1}}}^{2}Y_{{1,1}}
+{X_{{1}}}^{2}Y_{{1,2}}
+{Y_{{1,1}}}^{2}X_{{1}}
+{Y_{{1,1}}}^{2}Y_{{1,2}}
+{Y_{{1,2}}}^{2}X_{{1}}
+{Y_{{1,2}}}^{2}Y_{{1,1}}   \cr
&&+{Z_{{1,2,1}}}^{2}Z_{{1,1,1}}
+{Z_{{1,2,1}}}^{2}Z_{{1,2,2}}
+{Z_{{1,2,2}}}^{2}Z_{{1,1,1}}
+{Z_{{1,1,1}}}^{2}Z_{{1,2,2}} \cr
&&
+Y_{{1,2}}{Z_{{1,2,1}}}^{2}
+Y_{{1,2}}{Z_{{1,2,2}}}^{2}
+Y_{{1,2}}{Z_{{1,1,2}}}^{2}
+{Y_{{1,2}}}^{2}Z_{{1,2,1}} \cr
&&
+{Y_{{1,2}}}^{2}Z_{{1,2,2}}
+{Y_{{1,2}}}^{2}Z_{{1,1,1}}
+{Y_{{1,2}}}^{2}Z_{{1,1,2}}
+Z_{{1,1,2}}{Z_{{1,2,1}}}^{2}  \cr
&&
+Z_{{1,1,2}}{Z_{{1,2,2}}}^{2}
+Z_{{1,1,2}}{Z_{{1,1,1}}}^{2}
+Y_{{1,2}}{Z_{{1,1,1}}}^{2}
+ {Z_{{1,1,2}}}^{2}Z_{{1,2,1}} \cr
&&
+ {Z_{{1,1,2}}}^{2}Z_{{1,1,1}}
+ {Z_{{1,2,2}}}^{2}Z_{{1,2,1}}
+ Y_{{1,1}}{Z_{{1,2,2}}}^{2}
+ Y_{{1,1}}{Z_{{1,1,1}}}^{2} \cr
&&
+ Y_{{1,1}}{Z_{{1,1,2}}}^{2}
+{Z_{{1,1,2}}}^{2}Z_{{1,2,2}}
+{Y_{{1,1}}}^{2}Z_{{1,1,1}}
+{Y_{{1,1}}}^{2}Z_{{1,1,2}}   \cr
&&
+{Y_{{1,1}}}^{2}Z_{{1,2,2}}
+{Y_{{1,1}}}^{2}Z_{{1,2,1}}
+X_{{1}}{Z_{{1,1,2}}}^{2}
+X_{{1}}{Z_{{1,1,1}}}^{2}
+X_{{1}}{Z_{{1,2,2}}}^{2}   \cr
&&
+Y_{{1,1}}{Z_{{1,2,1}}}^{2}
+X_{{1}}{Z_{{1,2,1}}}^{2}
+{X_{{1}}}^{2}Z_{{1,1,2}}
+{X_{{1}}}^{2}Z_{{1,1,1}}
+{X_{{1}}}^{2}Z_{{1,2,2}}  \cr
&&
+{X_{{1}}}^{2}Z_{{1,2,1}}
+{Z_{{1,1,1}}}^{2}Z_{{1,2,1}}  \cr
&&
+2 \left( \right. \,Y_{{1,1}}Y_{{1,2}}Z_{{1,1,1}}
+ \,X_{{1}}Y_{{1,1}}Z_{{1,2,1}}
+ \,X_{{1}}Y_{{1,1}}Z_{{1,2,2}}
+ \,X_{{1}}Y_{{1,1}}Z_{{1,1,1}} \cr
&&
+ \,X_{{1}}Y_{{1,1}}Z_{{1,1,2}}
+ \,Y_{{1,1}}Y_{{1,2}}Z_{{1,1,2}}
+ \,Y_{{1,1}}Y_{{1,2}}Z_{{1,2,1}}
+ \,Y_{{1,1}}Y_{{1,2}}Z_{{1,2,2}}\cr
&&
+ \,X_{{1}}Y_{{1,2}}Z_{{1,2,2}}
+ \,X_{{1}}Y_{{1,2}}Z_{{1,1,1}}
+ \,X_{{1}}Y_{{1,2}}Z_{{1,1,2}}
+ \,X_{{1}}Y_{{1,2}}Z_{{1,2,1}}  \cr
&&
+ \,X_{{1}}Y_{{1,1}}Y_{{1,2}}
+ \,X_{{1}}Z_{{1,2,1}}Z_{{1,2,2}}
+ \,X_{{1}}Z_{{1,1,1}}Z_{{1,2,1}}\cr
&&
+ \,X_{{1}}Z_{{1,1,2}}Z_{{1,2,2}}
+ \,X_{{1}}Z_{{1,1,1}}Z_{{1,2,2}}
+ \,Z_{{1,1,2}}Z_{{1,2,1}}Z_{{1,2,2}} \cr
&&
+ \,Z_{{1,1,2}}Z_{{1,1,1}}Z_{{1,2,1}}
+ \,Z_{{1,1,2}}Z_{{1,1,1}}Z_{{1,2,2}}
+
\,Y_{{1,1}}Z_{{1,2,1}}Z_{{1,2,2}}  \cr
&&
+ \,Y_{{1,1}}Z_{{1,1,1}}Z_{{1,2,1}}
+ \,Y_{{1,1}}Z_{{1,1,1}}Z_{{1,2,2}}
+ \,Y_{{1,1}}Z_{{1,1,2}}Z_{{1,2,2}} \cr
&&
+ \,Y_{{1,1}}Z_{{1,1,2}}Z_{{1,2,1}}
+ \,Y_{{1,1}}Z_{{1,1,2}}Z_{{1,1,1}}
+ \,X_{{1}}Z_{{1,1,2}}Z_{{1,2,1}}  \cr
&&
+ \,X_{{1}}Z_{{1,1,2}}Z_{{1,1,1}}
+ \,Z_{{1,1,1}}Z_{{1,2,1}}Z_{{1,2,2}}
+ \,Y_{{1,2}}Z_{{1,1,1}}Z_{{1,2,1}}  \cr
&&
+ \,Y_{{1,2}}Z_{{1,2,1}}Z_{{1,2,2}}
+ \,Y_{{1,2}}Z_{{1,1,2}}Z_{{1,1,1}}
+ \,Y_{{1,2}}Z_{{1,1,2}}Z_{{1,2,1}}  \cr
&&
+ \,Y_{{1,2}}Z_{{1,1,2}}Z_{{1,2,2}}
+ \,Y_{{1,2}}Z_{{1,1,1}}Z_{{1,2,2}}
   \left.\right)\; \left.\right).
\nonumber \end{eqnarray}

\begin{sidewaystable}
{\footnotesize
\caption{
 Data of indeterminates relative to multi-indicial symmetric
functions.}\label{tab1} $\\\\$
\begin{center}
 \begin{tabular}{|l||l||llll||l||l|}
 \hline &&&$\left\{ {{\bf a}}_{ \left\{ 1\leq p \leq m,
\left[\mu\right]_{1\leq i\leq k}\right\} } \right\}$&&& $\longrightarrow$ Adding $m+1$
& $\longleftarrow$ $h$-projection \\
&&&&&&& on $\Lambda_{m,k} $\\
 \hline\hline &$D^0_{(m)}$& ${\bf a} _{1}$ &
${\bf a}_{2}$ & $\cdots$ &   ${\bf a}_{m}$ &- & ${\bf a}_{m+1}$\\
 &$D^1_{(m)}$&
${\bf a}_{1\mu}$ & ${\bf a}_{2\mu}$ &   $\dots$ & ${\bf a}_{m \mu}$ & $1\leq \mu \leq n$ &
${\bf a}_{(m+1) \mu}$
\\
&
$D^2_{(m)}$& ${\bf a}_{1{\mu}_{1}{\mu}_{2}}$ &
${\bf a}_{2{\mu}_{1}{\mu}_{2}}$ & $\dots$ & ${\bf a}_{m{\mu}_{1}{\mu}_{2}}$ &  $1\leq
{\mu}_1,{\mu}_2\leq n$ & ${\bf a}_{(m+1){\mu}_{1}{\mu}_{2}}$\cr
&.&. &.  & $\dots$ &. &. &.  \cr
&. &. &. & $\dots$ &. &. &. \cr
&. &. &. & $\dots$ &. &. &. \cr
& $D^k_{(m)}$&   ${\bf a}_{1{\mu}_{1}{\mu}_{2}\dots{\mu}_{k}}$
& ${\bf a}_{2{\mu}_{1}{\mu}_{2}\dots{\mu}_{k}}$ &   $\dots$ &
${\bf a}_{m{\mu}_{1}{\mu}_{2}\dots{\mu}_{k}}$ &   $1\leq {\mu}_1,{\mu}_2,\dots
,{\mu}_{k}\leq n$ &   ${\bf a}_{(m+1){\mu}_{1}{\mu}_{2}\dots{\mu}_{k}}$
\cr
\hline
$v$-Projection $\uparrow$&-&
$D^{(k)}_1$ &   $D^{(k)}_2$ &   $\dots$ &   $D^{(k)}_m$ &- &$D^{(k)}_{m+1}$\cr
on $\Lambda_{m,k}$&&&&& &&\cr
\hline
Adding $k+1\;\;$ $\downarrow$&
$D^{k+1}_{(m)}$ &
${\bf a}_{1{\mu}_{1}{\mu}_{2}\dots{\mu}_{k}{\mu}_{k+1}}$ &
${\bf a}_{2{\mu}_{1}{\mu}_{2}\dots{\mu}_{k}{\mu}_{k+1}}$ &   $\dots$ &
${\bf a}_{m{\mu}_{1}{\mu}_{2}\dots{\mu}_{k}{\mu}_{k+1}}$ &   $1\leq {\mu}_1,{\mu}_2,\dots ,$ &
${\bf a}_{{(m+1)}{\mu}_{1}{\mu}_{2}\dots{\mu}_{k}{\mu}_{k+1}}$
\cr
&&&&&& ${\mu}_{k},{\mu}_{k+1}\leq n$ & \cr
\hline
\end{tabular}\end{center}}
\end{sidewaystable}
\begin{sidewaystable}
{
{\caption{Commutative diagram of general inverse limits}\label{tab2} $\\\\$
\begin{center}
\begin{tabular}{|l||llllllllll|l|}
\hline
 &  &$\Lambda$  & &$\stackrel{H_{m,.}}{\cdots\longrightarrow}$&&$ \Lambda_{m,.}$       &&&&&  \\
\hline\hline
$\Lambda$  &   &$\Lambda$ & $\cdots >$&
$\Lambda_{m+p,.}$&$\stackrel{\phi_{(m+p,m)}}{\cdots\longrightarrow}$ &$ \Lambda_{m,.}$& &&
$\stackrel{h_m}{\longrightarrow}$&$\Lambda_{m-1,.}$&$\Lambda_{m-1,.}$\\
  &  & $\stackrel{\vdots}{\vee}$   & & & &$\vdots$ & & &  &$\vdots$ &   \\

 $V_{k}\stackrel{\vdots}{\downarrow}$ &  & $\Lambda_{.,(k+q)}$ &$\cdots\longrightarrow$&
  $\Lambda_{m+p, k+q}$&$\cdots\longrightarrow$&$\Lambda_{m, k+q}$ && &$\stackrel{h_{m,k+q}}{\longrightarrow}$
  & $\Lambda_{m-1,k+q}$
  & $\stackrel{\vdots}{\downarrow} V_{m-1,k}$\\

   & $\psi_{(k+q,k)}$ & $\stackrel{\vdots}{\downarrow}$ &
 &&$\psi_{m,(k+q,k)}$&$\stackrel{\vdots}{\downarrow}$&&&&
 $\stackrel{\vdots}{\downarrow}\psi_{m-1,(k+q,k)}$ &\\

 & &  &  &&&&&&&&\\

 $\Lambda_{.,k}$& &$\Lambda_{.,k}$&$\cdots\longrightarrow$
 & $\Lambda_{m+p, k}$&$\cdots\longrightarrow$&$\Lambda_{m, k}$
 &&&$\stackrel{h_{m,k}}{\longrightarrow}$&$\Lambda_{m-1, k}$& $\Lambda_{m-1, k} $\\

 & $v_k$ & $\downarrow$ & & &  &$v_{m,k}\downarrow$&  & &&$\downarrow v_{m-1,k}$ &  \\

 &  & $\Lambda_{.,k-1}$& $\cdots \longrightarrow$  & $\Lambda_{m+p, k-1}$
 &$\cdots\longrightarrow$&$\Lambda_{m, k-1}$&
 &&$\stackrel{h_{m,k-1}}{\longrightarrow}$&$\Lambda_{m-1, k-1}$ &  \\
\hline
 &  & $\Lambda_{.,k-1}$ && $\stackrel{\cdots\longrightarrow}{H_{m,k-1}}$ && $\Lambda_{m, k-1}$ &&&& &  \\
\hline
\end{tabular}
\end{center} }}
\end{sidewaystable}
\begin{sidewaystable}
{\scriptsize
{\caption{Matrix ${\mathcal A}$ generating the $\left[ m,k \right]$ Vandermonde
determinant.}\label{tab3}
\begin{center}$f(p,0)=p;\; 1\leq p \leq m;\;\;\;\;
f(p, \left[\mu\right]_t)\equiv f(p, \mu_1, \mu_2, \dots, \mu_t )
=mq_{t-1} + (q-1)n^{t} +\sum_{l=1}^{t}(\nu_l-1)n^{t-l}+\nu_q$,
$t\geq 1,\; 1\leq\nu_l\leq n,$\end{center}
\begin{center}$q_{k}=\frac{n^{k+1}-1}{n-1},\; n\neq 1;\;\; q_{k}=k+1,\; n=1.$
\end{center}
$\\\\$
\begin{center}
\begin{tabular}{|l|lllll|}\hline Row&&Column &
&Column  &Column
\\
&& $1\leq q\leq m$& &$f(q,\left[\mu\right]_t)$ & $f(m,n\dots
n)$\\
&& & & & $=mq_k$
\\
\hline $f(1,0)=1:$&&$ {\bf a}_{1}^{\delta^{\left[ m,0\right]}_{1}},
{\bf a}_{1}^{\delta^{\left[ m,0\right]}_{2}},\dots, {\bf
a}_{1}^{\delta^{\left[ m,0\right]}_{m}},$&$ {\bf
a}_{1}^{\delta^{\left[ m,1\right]}_{11}}, \dots,{\bf
a}_{1}^{\delta^{\left[ m,1\right]}_{1n}}, {\bf
a}_{1}^{\delta^{\left[ m,1\right]}_{21}}, \dots,{\bf
a}_{1}^{\delta^{\left[ m,1\right]}_{2n}},\dots,$&$ {\bf
a}_{1}^{\delta^{\left[ m,t \right]}_{q\nu_1\nu_2\dots\nu_t}},
\dots,$&${\bf a}_1^{\delta^{\left[m,k\right]}_{mnn\dots n}}$\\
&$\vdots$&&&&\\
$f(p,0)=p$:&& ${\bf a}_p^{\delta^{\left[ m,0\right]}_{1}},{\bf
a}_p^{\delta^{\left[ m,0\right]}_{2}},\dots,$&$\dots$ &${\bf
a}_{p}^{\delta^{\left[ m,t \right]}_{q\nu_1\nu_2\dots\nu_t}},
\dots,$&${\bf a}_{p}^{\delta^{\left[m,k\right]}_{mnn\dots n}}$\\
&$\vdots$&&&&\\
$f(m,0)=m$:&& $ {\bf a}_{m}^{\delta^{\left[ m,0\right]}_{1}}, {\bf
a}_{m}^{\delta^{\left[ m,0\right]}_{2}},\dots,$ &$\dots$& ${\bf
a}_{m}^{\delta^{\left[ m,t \right]}_{q\nu_1\nu_2\dots\nu_t}},
\dots,$&${\bf a}_{m}^{\delta^{m,k}_{mnn\dots n}}$\\
$f(m,1,0,\dots)=m+1$:&& ${\bf a}_{11}^{\delta^{\left[
m,0\right]}_{1}}, {\bf a}_{11}^{\delta^{\left[
m,0\right]}_{2}},\dots,$&$\dots$& ${\bf a}_{11}^{\delta^{\left[ m,t
\right]}_{q\nu_1\nu_2\dots\nu_t}},
\dots,$&${\bf a}_{11}^{\delta^{\left[m,k\right]}_{mnn\dots n}}$\\
&$\vdots$&&&&\\
$f(m,\mu,0,\dots)=m+\mu$:&&$ {\bf a}_{1\mu}^{\delta^{\left[
m,0\right]}_{1}}, {\bf a}_{1\mu}^{\delta^{\left[
m,0\right]}_{2}},\dots,$&$\dots$ &${\bf a}_{1\mu}^{\delta^{\left[
m,t \right]}_{q\nu_1\nu_2\dots\nu_t}},
\dots,$&${\bf a}_{1\mu}^{\delta^{\left[m,k\right]}_{mnn\dots n}}$\\
&$\vdots$&&&& \\
$f(p,\left[\mu\right]_{1})=(p-1)n+\mu+m $:&& $ {\bf
a}_{p\mu}^{\delta^{\left[ m,0\right]}_{1}}, {\bf
a}_{p\mu}^{\delta^{\left[ m,0\right]}_{2}},\dots,$&$\dots$& $ {\bf
a}_{p\mu}^{\delta^{\left[ m,t \right]}_{q\nu_1\nu_2\dots\nu_t}}
, \dots,$&${\bf a}_{p\mu}^{\delta^{\left[m,k\right]}_{mnn\dots n}}$\\
&$\vdots$&&&&\\
$f(p,\left[\mu \right]_2)$: &&${\bf a}_{p\mu_1\mu_2}^{\delta^{\left[
m,0\right]}_{1}},
 {\bf a}_{p\mu_1\mu_2}^{\delta^{\left[
m,0\right]}_{2}},\dots,$&$\dots$&
$ {\bf a}_{p\mu_1\mu_2}^{\delta^{\left[ m,t \right]}_{q\nu_1\nu_2\dots\nu_t}},
\dots,$&${\bf a}_{p\mu_1\mu_2}^{\delta^{\left[m,k\right]}_{mnn\dots n}}$\\
&$\vdots$&&&&\\
$f(p, \left[\mu\right]_d)$: && ${\bf
a}_{p\mu_{1}\mu_{2}\dots\mu_{d}}^{\delta^{\left[ m,0\right]}_{1}},
{\bf a}_{p\mu_{1}\mu_{2}\dots\mu_{d}}^{\delta^{\left[
m,0\right]}_{2}},\dots,$&$\dots$&$ {\bf
a}_{p\mu_{1}\mu_{2}\dots\mu_{d}}^{\delta^{\left[ m,t
\right]}_{q\nu_1\nu_2\dots\nu_t}},
\dots,$&${\bf a}_{p\mu_{1}\mu_{2}\dots\mu_{d}}^{\delta^{\left[m,k\right]}_{mnn\dots n}}$\\
&$\vdots$&&&&\\
$f(m,\left[\mu\right]_{k})$:&&
 ${\bf a}_{m\mu_{1}\mu_{2}\dots\mu_{k}}^{\delta^{\left[
m,0\right]}_{1}}, {\bf a}_{m\mu_{1}\mu_{2}\dots\mu_{k}}^{\delta^{\left[ m,0\right]}_{2}},\dots,$&\dots&$
{\bf a}_{m\mu_{1}\mu_{2}\dots\mu_{k}}^{\delta^{\left[ m,t \right]}_{q\nu_1\nu_2\dots\nu_t}},
\dots,$&${\bf a}_{m\mu_{1}\mu_{2}\dots\mu_{k}}^{\delta^{\left[m,k\right]}_{mnn\dots n}}$\\
&$\vdots$&&&&\\
$f(m,n\dots n)=mq_k$:&&${\bf a}_{mnn\dots n}^{\delta^{\left[
m,0\right]}_{1}}, {\bf a}_{mnn\dots n}^{\delta^{\left[
m,0\right]}_{2}},\dots,$&\dots&$ {\bf a}_{mnn\dots
n}^{\delta^{\left[ m,t \right]}_{q\nu_1\nu_2\dots\nu_t}},
\dots,$&${\bf a}_{mnn\dots n}^{\delta^{\left[m,k\right]}_{mnn\dots
n}}$\\\hline
\end{tabular} \end{center}
 }}
\end{sidewaystable}


\begin{thebibliography}{[SiW]}

  \bibitem{IGM}
   I.~G.~Mcdonald,
   {\it Symmetric Functions and Orthogonal Polynomials},
   {\it Dean Jacqueline B. Lewis Memorial Lectures}, University
   Lecture Series , Vol. 12,
   Rutgers University, American Mathematical Society, 1998.

  \bibitem{GR}
   I.~Gelfand, D.~Krob, A.~Lascoux, B.~Leclerc,
   V.~S.~Retakh and J.~Y.~Thibon,
   {\it Noncommutative Symmetric Functions},
   Adv. in Math. {\bf 112 } 2 (1995), 218--348.

  \bibitem{GR2}
   I.~Gelfand, S.~Gelfand, V.~Retakh and R.~Lee~Wilson,
   {\it  Quasideterminants},
   to appear in Adv. in Math. {\bf 193} 1  (2005), 56--141;
   e-print {\tt arXiv: math.QA/0208146}.

    \bibitem{mac}
   P.~A.~McMahon, {\it Combinatory Analysis}, Cambridge University Press
   1915, 1916; Chelsea reprint 1960.

   \bibitem{gessel}
   I. M. Gessel, {\it Enumerative Applications of Symmetric Functions},
  Actes 17$^e$ Sem. Lothar. Combin. (1987), 5--21.

  \bibitem{dalbec}
  J. Dalbec, {\it Multisymmetric functions}, Beitrage Algebra Geom. {\bf 40} 1 (1999),
   27--51.

 \bibitem{vaccar}
 F. Vaccarino, {\it The ring of multisymmetric functions},
 Ann. Inst. Fourier {\bf 55} 3 (2005), 717--731.

 \bibitem{fresc}
 M. Feschbach, {\it The mod 2 cohomology rings of the symmetric group and invariants},
 Topology {\bf 41} (2002), 57--84.

  \bibitem{L}
    S.~Lang,
   {\it Algebra },
   $2^{nd}$Ed., Addison Wesley Publishing
   $C^{ie}$ inc.,
   Yale Univ. New Haven, Connecticut, U.S.A., 1984.

  \bibitem{bh}
   J.~Ben~Geloun and M.~N.~Hounkonnou,
   {\it in progress}.

\end{thebibliography}
\end{document}